%% file: Calo_Romkes_Valseth_AVS_FE_LNCSE.tex
\documentclass[graybox]{svmult}

\usepackage{theorem}
\usepackage{amssymb,amsmath}
\usepackage{cite}
\usepackage{mathptmx}       
\usepackage{helvet}         
\usepackage{courier}        
\usepackage{type1cm}        
\usepackage{makeidx}         
\usepackage{graphicx}        
\usepackage{multicol}        
\usepackage[bottom]{footmisc}
\usepackage{subfigure}

\newcommand{\bfm}[1]         { \mathbf {#1} }
\newcommand{\DD}            { {\bfm{D}}   }              
\newcommand{\bb}             { {\bfm{b}}   }              
\newcommand{\Nabla}       { \boldsymbol{\nabla} }   
\newcommand{\xx}             { \bfm{x}     }             
\newcommand{\qq}             { \bfm{q}     }             
\newcommand{\vn}             { \bfm{n}     }             

\newcommand{\GD}          { \Gamma_{D} }            
\newcommand{\GN}          { \Gamma_{N} }            
\newcommand{\Ph}            { \mathcal{P}_h }
\newcommand{\Kep}           { K_m \in \Ph}
\newcommand{\dKm}           { \partial K_m }

\newcommand{\FF}             { \bfm{F}_m   }      
\newcommand{\KK}           { \hat{K} }

\newcommand{\SCZO}            { C^0(\Omega) }                                
\newcommand{\SCZOvec}            { [C^0(\Omega)]^2 }                                
\newcommand{\SLTO}            { L^2(\Omega) }                                
\newcommand{\SLTOvec}            { [L^2(\Omega)]^2 }                                
\newcommand{\SHOO}            { H^1(\Omega) }                              
\newcommand{\SHOP}            { H^1(\Ph) }                              
\newcommand{\VV}            { {V(\Ph)} }                              
\newcommand{\VVK}            { {V(K_m)} }                              
\newcommand{\UU}            { U(\Omega) }                              
\newcommand{\VVh}            { V^*(\Ph) }                              
\newcommand{\UUh}            { U^h(\Omega) }                              
\newcommand{\SHOK}            { H^1(K_m) }                              
\newcommand{\SLTK}            { L^2(K_m) }                              
\newcommand{\SLTKvec}            { [L^2(K_m)]^2 }                              
\newcommand{\SHdivO}            { H(\text{div},\Omega) }                              
\newcommand{\SHdivK}            { H(\text{div},K_m) }                              
\newcommand{\SHMHGN}            { H^{-1/2}(\GN) }                              
\newcommand{\SHHdK}            { H^{1/2}(\dKm) }                              
\newcommand{\SHmHdK}            { H^{-1/2}(\dKm) }                              
\newcommand{\SLOinf}            { L^\infty(\Omega)}                               
\newcommand{\norm}[2]         { \| {#1} \|_{#2} }                      
\newcommand{\dx}              { \; {\rm d} \bfm{x}   }                                
\newcommand{\dss}             { \, {\rm d} s   }                                         
\newcommand{\summa}[2]        { \overset{#2}{\underset{#1}{\sum}} } 
\newcommand{\supp}[1]         { \underset{#1}{\sup} \, }                        
\newcommand{\wwwh}              { \bfm{w^*} }                      
\newcommand{\vvh}              { v^* }                      
\newcommand{\vvi}              { \tilde{e}^i }                      
\newcommand{\wwwi}              { \mathbf{\tilde{E}^i}}                     
\newcommand{\vvj}              { \tilde{e}^j_x }                      
\newcommand{\wwwj}              { \mathbf{\tilde{E}^j_x}}                     
\newcommand{\vvk}              { \tilde{e}^k_y }                      
\newcommand{\wwwk}              { \mathbf{\tilde{E}^k_y}}                     
\newcommand{\vvih}              { \tilde{e}^i_h }                      
\newcommand{\wwwih}              { \mathbf{\tilde{E}^i_h}}                     
\newcommand{\vvjh}              { \tilde{e}^{j}_{x_h} }                      
\newcommand{\wwwjh}              { \mathbf{\tilde{E}^{j}_{x_h}}}                     
\newcommand{\vvkh}              { \tilde{e}^k_{y_h} }                      
\newcommand{\wwwkh}              { \mathbf{\tilde{E}^k_{y_h}}}                     
\newcommand{\zzz}              { \bfm{z} }                      
\newcommand{\www}             { \bfm{w} }                      
\newcommand{\zz}              { \bfm{z} }                      
\newcommand{\nn}              { \bfm{n} }                      
\newcommand{\vfr}              { {\pmb \theta} }                      
\newcommand{\fr}              { { \varphi} }                      
\newcommand{\isdef}           { \overset{\text{def}}{=} }     
\newcommand{\ds}              { \displaystyle }               

\theoremstyle{plain}
\theoremheaderfont{\bfseries \upshape}

\newtheorem{lem}{Lemma}[section]
\newtheorem{rem}{Remark}[section]

\makeindex             

\begin{document}

\title*{Automatic Variationally Stable Analysis for FE Computations: An Introduction }
\author{Victor M. Calo, Albert Romkes, and Eirik Valseth}
\institute{Victor M. Calo \at Applied Geology Department, Curtin University, GPO Box U1987, Perth, Western Australia 6845, \email{victor.calo@curtin.edu.au}
\and Albert Romkes   (Corresponding Author) \at Department of Mechanical Engineering, South Dakota School of Mines \& Technology, 501 E. St. Joseph Street, Rapid City, SD 57701, USA, \email{Albert.Romkes@sdsmt.edu}
\and Eirik  Valseth \at Department of Mechanical Engineering, South Dakota School of Mines \& Technology, 501 E. St. Joseph Street, Rapid City, SD 57701, USA, \email{Eirik.Valseth@mines.sdsmt.edu}}%

%
\maketitle
\abstract*{We introduce an automatic variationally stable analysis (AVS)  for  finite element (FE) computations of scalar-valued convection-diffusion equations with non-constant and highly oscillatory coefficients.
In the spirit of least squares FE methods~\cite{bochevLeastSquares}, the AVS-FE method recasts the governing second order partial differential equation (PDE) into a system of first-order PDEs. However, in the subsequent derivation of the equivalent weak formulation,  a Petrov-Galerkin technique is applied by using different regularities for the trial and test function spaces. We use standard FE approximation spaces  for the trial spaces, which are  $C^0$, and broken  Hilbert spaces for the test functions. Thus, we seek to compute pointwise continuous solutions for both the primal variable and its flux (as in least squares FE methods), while the test functions are piecewise discontinuous.
To ensure the numerical stability of the subsequent FE discretizations,  we apply the  philosophy of the discontinuous Petrov-Galerkin (DPG) method by Demkowicz and Gopalakrishnan~\cite{Demkowicz4, Demkowicz1, Demkowicz2, Demkowicz3, Demkowicz5, Demkowicz6},  by invoking test functions that lead to unconditionally stable numerical systems (if the kernel of the underlying differential operator is trivial).
In the AVS-FE method, the discontinuous test functions are ascertained per the DPG approach from local, decoupled, and well-posed variational problems, which lead to best approximation properties in terms of the energy norm. 
We present various 2D numerical verifications, including convection-diffusion problems with highly oscillatory coefficients and extremely high Peclet numbers, up to $O(10^{9})$. These show the unconditional stability without the need for any upwind schemes nor any other artificial numerical stabilization. The results are not highly diffused  for convection-dominated problems nor show any strong oscillations, but adequately capture and indicate the presence of boundary layers, even for very coarse meshes and low polynomial degrees of approximation, $p$. Remarkably, we can compute the test functions by  using the same $p$ level as  the trial functions without significantly impacting the numerical accuracy or asymptotic convergence of the numerical results.  In addition, the AVS method delivers high numerical accuracy for the computed flux.  Importantly, the AVS methodology delivers optimal asymptotic error convergence rates of order $p+1$ and $p$ are obtained in the $L^2$ and $H^1$ norms for the primal variable. Our experience indicates that for convection-dominated problems we often observe a convergence rate of $p+1$   for the $L^2$ norm of the flux variable.}
\abstract{We introduce an automatic variationally stable analysis (AVS)  for  finite element (FE) computations of scalar-valued convection-diffusion equations with non-constant and highly oscillatory coefficients.
In the spirit of least squares FE methods~\cite{bochevLeastSquares}, the AVS-FE method recasts the governing second order partial differential equation (PDE) into a system of first-order PDEs. However, in the subsequent derivation of the equivalent weak formulation,  a Petrov-Galerkin technique is applied by using different regularities for the trial and test function spaces. We use standard FE approximation spaces  for the trial spaces, which are  $C^0$, and broken  Hilbert spaces for the test functions. Thus, we seek to compute pointwise continuous solutions for both the primal variable and its flux (as in least squares FE methods), while the test functions are piecewise discontinuous.
To ensure the numerical stability of the subsequent FE discretizations,  we apply the  philosophy of the discontinuous Petrov-Galerkin (DPG) method by Demkowicz and Gopalakrishnan~\cite{Demkowicz4, Demkowicz1, Demkowicz2, Demkowicz3, Demkowicz5, Demkowicz6},  by invoking test functions that lead to unconditionally stable numerical systems (if the kernel of the underlying differential operator is trivial).
In the AVS-FE method, the discontinuous test functions are ascertained per the DPG approach from local, decoupled, and well-posed variational problems, which lead to best approximation properties in terms of the energy norm. 
We present various 2D numerical verifications, including convection-diffusion problems with highly oscillatory coefficients and extremely high Peclet numbers, up to $O(10^{9})$. These show the unconditional stability without the need for any upwind schemes nor any other artificial numerical stabilization. The results are not highly diffused  for convection-dominated problems nor show any strong oscillations, but adequately capture and indicate the presence of boundary layers, even for very coarse meshes and low polynomial degrees of approximation, $p$. Remarkably, we can compute the test functions by  using the same $p$ level as  the trial functions without significantly impacting the numerical accuracy or asymptotic convergence of the numerical results.  In addition, the AVS method delivers high numerical accuracy for the computed flux.  Importantly, the AVS methodology delivers optimal asymptotic error convergence rates of order $p+1$ and $p$ are obtained in the $L^2$ and $H^1$ norms for the primal variable. Our experience indicates that for convection-dominated problems we often observe a convergence rate of $p+1$   for the $L^2$ norm of the flux variable.}
\section{Introduction }
\label{sec:introduction}
Singularly perturbed problems are ubiquitous in many engineering applications. We seek to develop a framework to tackle this large class of problems in a constructive manner. We start with a common model problem, that is, the convection-diffusion  problem which is relevant to many engineering applications where transport mechanisms play a significant role, e.g., subsurface flow through porous media, dynamics of viscous  flow, convective transfer of heat, drug delivery, turbulence modeling, etc.  In this paper, we focus on the stationary version of the scalar-valued convection-diffusion equation and therefore limit our consideration to  solutions which only depend on spatial variables and not the temporal variable. To date, the numerical analysis of even the stationary problem  poses  significant challenges due to the presence of the convection term, which dominates the diffusion processes. Classical FE methodologies, such as the Bubnov-Galerkin FE method~\cite{Bubnov1913,hughes2012finite,oden2012introduction,reddy1993introduction}, mixed FE methods~\cite{BrezziMixed,RaviartThomas}, and Petrov-Galerkin method~\cite{petrov1940application}, struggle in their numerical analysis due to the numerical instability introduced by the convection term. The corresponding discrete systems of equations can be ill-posed (i.e., a discrete solution does not exist) or lead to either spurious solutions or solutions with severe oscillations. These generally do not tend  to attenuate with continued mesh refinements and/or enrichments until the boundary layers are resolved, which in many applications  is prohibitively expensive. The least squares FE methods (LSFEMs)~\cite{bochevLeastSquares}, the $k$-version of the FE method by Surana~et~al.~\cite{SuranaIJES2003,Surana2010,Ahmadi2009,Maduri2009},  and the DPG method by Demkowicz and Gopalakrishnan~\cite{Demkowicz4, Demkowicz1, Demkowicz2, Demkowicz3, Demkowicz5, Demkowicz6} resolve the numerical instability issues by choosing test/weight functions that lead to unconditionally stable systems of equations governing the FE discretizations.  However, in the case of LSFEM and the $k$-version FE method, the numerical solutions, while stable, can be overly diffusive, particularly for coarse mesh partitions, and therefore fail to indicate the presence and/or location of any sharp boundary layers or other local solution features. As a result, the corresponding adaptive mesh strategies can be ineffective in the presence of strong convection and require overly refined mesh partitions with large numbers of degrees of freedom to resolve boundary layers or other local phenomena. Contrarily, the DPG method does not  suffer from overly diffused solutions but also requires edge fluxes and traces (referred to as numerical fluxes and traces). The number of degrees of freedom, once you statically condense the degrees of freedom internal to each element, is similar to the count of the method we propose herein. In addition, although the DPG method provides unconditionally stable FE discretizations, the stabilization is problem-dependent. To ensure the numerical stability and asymptotic convergence of the FE process, the numerical fluxes and traces  have to be numerically stabilized through multiplication  by mesh dependent terms.
This stabilization  is akin to upwind-schemes used in other FE methodologies  and depends highly on the form/nature of the diffusion and convection coefficients. It  is therefore problem-dependent.

Another technique which enlarges the approximation, introduced in~\cite{Calo2016}, extends the use of the generalized multiscale finite elements to stabilize the advection-diffusion model problem.
Alternatively, stabilized finite element methods do not add extra degrees of freedom to the global system, but require problem specific modifications of the stabilization parameter.  The original stabilization technique is the streamlined-upwind Petrov-Galerkin (SUPG) stabilization, introduced by Brooks and Hughes~\cite{Brooks1982} for the Navier-Stokes system. Using the analytical framework proposed by Hughes~\cite{Hughes1995}, we can interpret many  stabilization methods as residual-based modifications of the discrete weak forms where a locally scaled  differential operator acts on the test function to weight the residual of each trial function. The multiscale interpretation of the stabilization process was illuminating and opened many application opportunities~\cite{Hughes1996,Hughes1998,Hughes2004}, but did not simplify the design process of the stabilization technique. Effectively, this design process is arduous, and problem specific. Among the many successful stabilized methods we cite several that were applied to the transport and Navier-Stokes equations~\cite{Hauke1994,Behr1993,Codina1998,Codina2000,Franca1992,Franca1992a,Franca1993,Hughes1986, Hughes1986a, Hughes1987,Hughes1989,Jansen1999,Juanes2005,Shakib1991,Shakib1991a}.

In this manuscript, we introduce the automatic variationally stable  (AVS) analysis  for FE computations of the convection-diffusion equation in which the diffusion and convection coefficients can be highly oscillatory. The  method is essentially a hybrid of the LSFEM, Petrov-Galerkin, and the DPG methods by employing the strength and benefits of each approach separately, leading to a FE process that is unconditionally stable and produces numerical solutions that are not overly diffusive, even for coarse FE mesh partitions and low polynomial degrees of approximation. There is no need for the determination of any mesh- and problem-dependent stabilization parameters to warrant unconditionally stable numerical schemes nor overly refined/enriched initial FE mesh partitions to ascertain the presence and location of any boundary layers or local  phenomena.

Firstly, we follow  mixed FE approaches, by introducing the fluxes  as  auxiliary variables and thereby recast the second order, scalar-valued convection, diffusion problem into a first order vector-valued PDE.  We subsequently apply the Petrov-Galerkin philosophy in the derivation of the equivalent integral formulation (i.e. the weak form) of the established vector-valued PDE by allowing a different regularity for the trial and  test spaces. The FE discretization of the weak form is then applied such that the base variable and the fluxes are classical global $C^0$ functions. However, we apply  broken (i.e., discontinuous) Hilbert spaces for the test functions in an effort to allow a maximum flexibility in choosing  test functions that lead to unconditionally stable FE processes. To do so, we invoke the philosophy of the DPG method in the FE discretization of the weak form   by constructing a test function for every $C^0$ trial function, which is a solution to decoupled element-wise local variational problems; called 'test problems.' Conforming to the DPG philosophy, the test problems employ bilinear forms which define a local inner product on each element. In the AVS-FE method, we apply a local $H^1$ inner product as the bilinear form in the test problems. As in  DPG, the resulting test functions lead to unconditionally stable systems of equations governing the FE approximation of the problem. In addition, the specific choice for local $H^1$ inner products in the test problems, appears to result in FE approximations that are not overly diffusive; even for convection-dominated problems with Peclet numbers of order  $10^{9}$. Remarkably, the numerical solutions we obtain for the flux variables with the AVS-FE method are highly accurate.

Our choice for  $C^0$ trial functions is motivated by the fact that it enables us to enforce the continuity of  all variables strongly and in a straightforward manner. This is of particular benefit for the analysis of the fluxes  in the presence of highly oscillatory diffusion coefficients. Moreover, it negates the need  to introduce numerical  (edge) fluxes and traces as auxiliary variables and thereby  reduces the computational cost and removes the need for any problem-dependent numerical stabilization of such variables. A key benefit of this functional choice is that legacy software for pre- and post-processing the data for the simulations can be directly used to prepare and analyze the data required and produced by  AVS-FE. Importantly, our simulations rely on continuous discretizations which facilitate solution interpretations and analyses from an engineering point of view.   That is,  from the user point of view, they are standard finite element solutions where all variables are continuous, simplifying the adoption of the technique by the end-user community.

As in the DPG method, we establish a best approximation property in terms of the energy norm that is induced by the bilinear form of the integral formulation of the AVS-FE method and obtain optimal asymptotic convergence rates in $L^2$ and $H^1$ for the base  variable and in $L^2$ for the flux variables.

In the following, we present the derivation of the AVS-FE weak formulation for the convection-diffusion problem in Section~\ref{sec:integral_form} and its subsequent FE discretization in Section~\ref{sec:discrete_form}, and various two-dimensional verifications in Section~\ref{sec:numerical_verifications}. Concluding remarks and future efforts are discussed in Section~\ref{sec:conclusions}.

\section{Derivation of Integral Statement and FE Discretization}
\label{sec:CDPG_method}
%
%
%
%
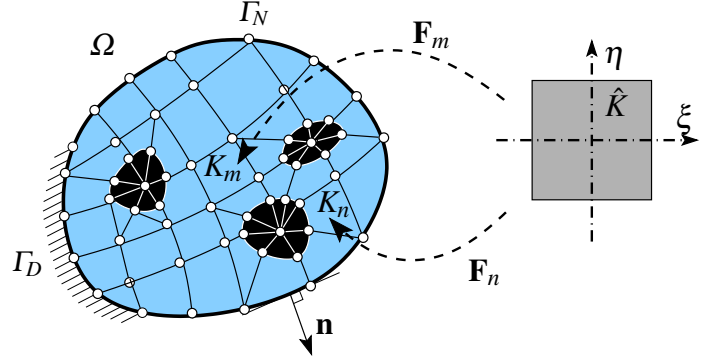
\begin{figure}[t]
\centering
\input{figures/fig-model-with-partition.pstex_t}
\caption{The model problem.}
\label{fig:model_problem}
\end{figure}
Let $\Omega \subset \mathbb{R}^2$ be an open bounded domain (see Figure~\ref{fig:model_problem}) with Lipschitz boundary $\partial \Omega$ and outward unit normal vector $\vn$. The boundary $\partial \Omega$ consists of  open subsections $\GD,\GN\subset\partial\Omega$, such that $\GD\cap\GN=\emptyset$ and $\partial \Omega = \overline{\GD\cup\GN}$. For our model problem, we  consider the following convection-diffusion equation in $\Omega$ with homogeneous Dirichlet boundary conditions applied on $\GD$ and (possibly) non-homogeneous Neumann boundary conditions on $\GN$:

\begin{equation} \label{eq:conv_diff_BVP}
\boxed{
\begin{array}{l}
\text{Find }  u  \text{ such that:}    
\\[0.05in] 
\qquad 
\begin{array}{rcl}
  -\Nabla \cdot(\DD \Nabla u) \, + \,
 \bb \cdot \Nabla u  & = & f, \quad \text{ in } \, \Omega, 
 \\[0.05in]
 \qquad u &  = & 0, \quad \text{ on } \, \GD , 
 \\
 \qquad \DD\Nabla u \cdot \vn & = & g, \quad \text{ on } \,  \GN,
 \end{array}
 \end{array}
}
\end{equation}
where $\DD$ denotes the second order diffusion tensor, with symmetric, bounded, and positive definite coefficients $D_{ij}\in\SLOinf$; $\bb\in[\SLTO]^2$ the convection coefficient; $f\in\SLTO$ the source function; and $g\in\SHMHGN$ the Neumann boundary data. We consider the scenario in which the diffusion coefficients $D_{ij}$ can be highly heterogeneous and therefore can change many orders in magnitude over small length scales throughout $\Omega$ (e.g., in Figure~\ref{fig:model_problem}, the differently colored subdomains represent areas with different values of the diffusions coefficients).

In this work, we seek to derive a DPG weak formulation of~\eqref{eq:conv_diff_BVP} by using a  regular partition $\Ph$ of $\Omega$ into  open subdomains, or elements, $K_m$ (see Figure~\ref{fig:model_problem}), with diameters $h_m$, such that :
\begin{equation}
\notag
\label{eq:domain}
  \Omega = \text{int} ( \bigcup_{\Kep} \overline{K_m} ).
\end{equation}
Any such partition $\Ph$ of $\Omega$ is applied  such that any discontinuities in the diffusion coefficient $D_{ij}$ or convection coefficient $\bb$ are restricted to the boundaries $\dKm$ of the elements $\Kep$ (see Figure~\ref{fig:model_problem}). That is, we assume our mesh fully resolves these spatial features, while it may not resolve the induced internal layers.

We  apply a mixed FE methodology and introduce the flux $\qq =\{q_x,q_y\}^T= \DD\Nabla u$ as an auxiliary variable, then, accordingly, $\qq\in\SHdivO$ and \eqref{eq:conv_diff_BVP} can be recast equivalently as a first-order system of PDEs, where the regularity of $u$ can be relaxed to be in $\SHOO$:

\begin{equation} \label{eq:conv_diff_BVP_first_order}
\boxed{
\begin{array}{l}
\text{Find }  (u,\qq) \in \SHOO\times\SHdivO \text{ such that:}    
\\[0.05in] 
\qquad 
\begin{array}{rcl}
  \qq - \DD\Nabla u & =  & 0, \quad \text{ in } \, \Omega, 
  \\
  -\Nabla \cdot \qq \, + \,
 \bb \cdot \Nabla u  & = & f, \quad \text{ in } \, \Omega, 
 \\[0.025in]
 \qquad u &  = & 0, \quad \text{ on } \, \GD , 
 \\
 \qquad \qq \cdot \vn & = & g, \quad \text{ on } \,  \GN. 
 \end{array}
 \end{array}
}
\end{equation}
\subsection{Derivation of Integral Formulation}
\label{sec:integral_form}
To start the derivation of the DPG formulation of~\eqref{eq:conv_diff_BVP_first_order}, we enforce the PDE weakly on each element $\Kep$, i.e., we seek the restrictions $u_m$ and $\qq_m$ of $u$ and $\qq$ to each $K_m$, such that:
\begin{equation} \label{eq:weak_BVP_L2}
\begin{array}{c}
\ds \int_{K_m} \biggl\{\left[ \qq_m - \DD\Nabla u_m\right] \cdot \www_m \, 
+  \, \left[  -\Nabla \cdot \qq_m \, + \,
 \bb \cdot \Nabla u_m \right] \, v_m \biggr\} \dx =   \int_{K_m} f\,v_m \dx, 
 \\[0.15in]
  \hfill \forall (v_m,\www_m) \in \SLTK\times\SLTKvec.
  \end{array}
\end{equation}
By repeating this process for all $\Kep$ and summing the resulting local integral formulations, we get: 
\begin{equation} \label{eq:weak_BVP_L2_sum}
\begin{array}{l}
\text{Find } \; (u,\qq) \in \SHOO\times\SHdivO:  
\\[0.1in] 
\ds \summa{\Kep}{}\int_{K_m} \biggl\{ \left[ \qq_m - \DD\Nabla u_m\right] \cdot \www_m \, 
+  \, \left[  -\Nabla \cdot \qq_m \, + \,
 \bb \cdot \Nabla u_m \right] \, v_m \biggr\} \dx  
 \\[0.125in]
 \hspace{0.75in} \ds =   \summa{\Kep}{}\int_{K_m} f\,v_m \dx, \qquad \forall (v,\www) \in \SLTO\times\SLTOvec
 \end{array}
\end{equation}
Next, we apply Green's identity to the $(\Nabla \cdot \qq_m) \,v_m$ terms, which demands that we increase the  regularity of each $v$ to be in $H^1$ locally for every $\Kep$. i.e.,
\begin{equation} \label{eq:weak_BVP_sum_1}
\begin{array}{l}
\text{Find } \; (u,\qq) \in \SHOO\times\SHdivO:  
\\[0.1in] 
\ds \summa{\Kep}{}\biggl\{ \int_{K_m}\biggl[  \, \left( \qq_m - \DD\Nabla u_m\right) \cdot \www_m \, 
+  \, \qq_m \cdot \Nabla v_m \, + \,
 (\bb \cdot \Nabla u_m) \, v_m  \biggr] \dx\biggr.
 \\[0.1in]
 \hspace{1in} \biggl. \ds - \oint_{\dKm}  \gamma^m_\nn(\qq_m) \, \gamma^m_0(v_m) \, \dss \biggr\}=   \summa{\Kep}{}\int_{K_m} f\,v_m \dx, 
 \\[0.2in]
 \hfill \forall (v,\www) \in \SHOP\times\SLTOvec
 \end{array}
\end{equation}
where the broken $H^1$ Hilbert space on $\Ph$ is defined as follows:
\begin{equation}
\label{eq:broken_h1_space}
\SHOP \isdef \biggl\{ v\in\SLTO: \quad v_m \in \SHOK, \; \forall \Kep\biggr\},
\end{equation}
and $\gamma^m_0: \SHOK: \longrightarrow \SHHdK$ and $\gamma^m_\nn:\SHdivK \longrightarrow \SHmHdK$ denote the trace and normal trace operators (e.g., see~\cite{Girault1986}) on $K_m$; and $\nn_m$ is the outward unit normal vector to the element boundary $\dKm$ of $K_m$. Strictly speaking, the edge integral on $\dKm$ in~\eqref{eq:weak_BVP_sum_1} is to be interpreted as the duality pairing in $\SHmHdK\times \SHHdK$ of $ \gamma^m_\nn(\qq_m)$ and $ \gamma^m_0(v_m)$, but we apply the engineering notation here by using an integral representation. 

Now, by decomposing each edge term in~\eqref{eq:weak_BVP_sum_1} into a sum of  several terms, i.e., one term concerning the portion of the edge $\dKm$ that intersects with neighboring elements and possibly one or two additional terms concerning the portion of $\dKm$ that intersects with $\GD$ or $\GN$, we can rewrite~\eqref{eq:weak_BVP_sum_1} as follows:
\begin{equation} \notag
\begin{array}{l}
\text{Find } \; (u,\qq) \in \SHOO\times\SHdivO:  
\\[0.1in] 
\ds \summa{\Kep}{}\biggl\{ \int_{K_m}\biggl[  \, \left( \qq_m - \DD\Nabla u_m\right) \cdot \www_m \, 
+  \, \qq_m \cdot \Nabla v_m \, + \,
 (\bb \cdot \Nabla u_m) \, v_m  \biggr] \dx\biggr.
 \\[0.1in]
 \hspace{0.5in} \ds 
 - \int_{\dKm\setminus \overline{\GD\cup\GN}}  \gamma^m_\nn(\qq_m) \, \gamma^m_0(v_m) \, \dss 
 - \int_{\dKm\cap\GD}  \gamma^m_\nn(\qq_m) \, \gamma^m_0(v_m) \, \dss
 \\[0.15in]
 \hspace{1in} \biggl. \ds 
 - \int_{\dKm\cap\GN}  \gamma^m_\nn(\qq_m) \, \gamma^m_0(v_m) \, \dss\biggr\}
 =   \summa{\Kep}{}\int_{K_m} f\,v_m \dx, 
 \\[0.2in]
 \hfill \forall (v,\www) \in \SHOP\times\SLTOvec
 \end{array}
\end{equation}
By subsequently enforcing the  Neumann boundary condition on the normal trace of $\qq$ as well as  constraining the traces of the test function $v_m$ on the Dirichlet boundary (since we apply the Dirichlet condition on $u$ strongly), we arrive at the final variational statement:
\begin{equation} \label{eq:weak_BVP_sum_2}
\begin{array}{l}
\text{Find } \; (u,\qq) \in \UU : 
\\[0.1in] 
\ds \summa{\Kep}{}\biggl\{ \int_{K_m}\biggl[  \, \left( \qq_m - \DD\Nabla u_m\right) \cdot \www_m \, 
+  \, \qq_m \cdot \Nabla v_m \, + \,
 (\bb \cdot \Nabla u_m) \, v_m  \biggr] \dx\biggr.
 \\[0.1in]
 \hspace{0.05in} \ds 
 - \int_{\dKm\setminus \overline{\GD\cup\GN}}  \gamma^m_\nn(\qq_m) \, \gamma^m_0(v_m) \, \dss \biggr\}
 =  \summa{\Kep}{}\left\{ \int_{K_m} f\,v_m \dx \, + \, \int_{\dKm\cap \GN} g \, \gamma^m_0(v_m) \dss \right\}, 
 \\[0.2in]
 \hfill \forall (v,\www) \in \VV
 \end{array}
\end{equation}
where the  trial and test function spaces, $\UU$ and $\VV$, are defined as follows:
\begin{equation}
\label{eq:function_spaces}
\begin{array}{c}
\UU \isdef \biggl\{ (u,\qq)\in \SHOO\times\SHdivO: \; \gamma_0^m(u_m)_{|\dKm\cap\GD} =0, \;  \forall\Kep\biggr\},
\\[0.15in]
\VV \isdef \biggl\{ (v,\www)\in \SHOP\times\SLTOvec: \,  \gamma_0^m(v_m)_{|\dKm\cap\GD} =0, \; \forall\Kep\biggr\},
\end{array}
\end{equation}
with  norms $\norm{\cdot}{\UU}:  \UU \!\! \longrightarrow \!\! [0,\infty)$ and $\norm{\cdot}{\VV}: \VV\! \! \longrightarrow\! \! [0,\infty)$ defined as:
\begin{equation}
\label{eq:broken_norms}
\begin{array}{l}
\ds \norm{(u,\qq)}{\UU} \isdef \sqrt{\int_{\Omega} \biggl[ \Nabla u \cdot \Nabla u + u^2   + (\Nabla \cdot \qq)^2+\qq \cdot \qq\biggr] \dx }.
\\[0.2in]
\ds   \norm{(v,\www)}{\VV} \isdef \sqrt{\summa{\Kep}{}\int_{K_m} \biggl[  h_m^2 \Nabla v_m \cdot \Nabla v_m + v_m^2   + \www_m \cdot \www_m\biggr] \dx }.
 \end{array}
\end{equation}
By introducing the bilinear form, $B:\UU\times\VV\longrightarrow \mathbb{R}$, and linear functional, $F:\VV\longrightarrow \mathbb{R}$, i.e.,
\begin{equation} \label{eq:B_and_F}
\begin{array}{c}
B((u,\qq);(v.\www)) \isdef
\ds \summa{\Kep}{}\biggl\{ \int_{K_m}\biggl[  \, \left( \qq_m - \DD\Nabla u_m\right) \cdot \www_m \, 
+  \, \qq_m \cdot \Nabla v_m \, + \,
 (\bb \cdot \Nabla u_m) \, v_m  \biggr] \dx\biggr.
 \\[0.1in]
 \hfill \ds 
 - \int_{\dKm\setminus \overline{\GD\cup\GN}}  \gamma^m_\nn(\qq_m) \, \gamma^m_0(v_m) \, \dss \biggr\},
 \\[0.15in]
 F((v,\www)) \isdef   \ds \summa{\Kep}{}\left\{\int_{K_m} f\,v_m \dx + \int_{\dKm\cap\GN} g\, \gamma^m_0(v_m) \dss\right\}, 
 \end{array}
\end{equation}
we can rewrite the  weak formulation~\eqref{eq:weak_BVP_sum_2} in compact form as follows:
\begin{equation} \label{eq:weak_form}
\boxed{
\begin{array}{ll}
\text{Find } (u,\qq) \in \UU & \hspace{-0.05in} \text{ such that:}
\\[0.05in]
 &  \quad B((u,\qq);(v,\www)) = F((v,\www)), \quad \forall (v,\www)\in \VV. 
 \end{array}}
\end{equation}
\begin{lem}
\label{lem:well_posed_weak_form}
Let $f\in(\SHOP)'$ and $g\in\SHMHGN$. Then there exists a unique solution $(u,\qq)\in \UU$ of the  weak formulation~\eqref{eq:weak_form}.
\end{lem}
We refer to~\cite{CaloRomkes2018} for a   proof of this lemma~\qed
\\

Now, \eqref{eq:weak_form} essentially represents a DPG formulation~\cite{Demkowicz4, Demkowicz1, Demkowicz2, Demkowicz3, Demkowicz5, Demkowicz6,NIEMI20132096}, as the spaces $\UU$ and $\VV$ have different regularities. However, it differs significantly  by using (weakly) globally continuous trial spaces. Currently existing DPG methods require weak enforcement of continuity conditions across inter-element edges by introducing numerical traces and fluxes as auxiliary variables. Thus, by  employing  trial spaces in which continuity of the primal variable and the normal fluxes is inherent (weakly), we attempt to keep the formulation, from the point of view of the user, as close as possible to a standard FE discretization.
 Lastly, the discrete description of the solution behaves like standard finite element discretizations, which will accelerate the adoption of this discretization technique by practitioners and paves the way to extend it to solutions with higher order global continuity, such as the ones produced by isogeometric analysis~\cite{hughes2005isogeometric,cottrell2009isogeometric,cottrell2006isogeometric,bazilevs2010isogeometric1,bazilevs2006isogeometric,cottrell2007studies,cortes2015performance,collier2014computational,puzyrev2017dispersion,bazilevs2009patient,bazilevs20113d,bazilevs2010large,hsu2012fluid,bazilevs2007variational,bazilevs2008isogeometric,gomez2008isogeometric,elguedj2008b,bazilevs2010isogeometric2,bazilevs2007weak,calo2008multiphysics,chang2012isogeometric,duddu2012finite,espath2016energy},  to show just a few of the relevant applications of this powerful simulation technique.

\subsection{AVS-FE Discretization}
\label{sec:discrete_form}
We now seek numerical approximations $(u^h,\qq^h)$ of solutions $(u,\qq)$ of  the weak form~\eqref{eq:weak_form} by using classical \emph{globally continuous}, $\SCZO$, trial functions for $(u^h,\qq^h)$.  However, the discontinuous topology is maintained for the space of test functions, as this allows the maximum flexibility in  constructing test functions that lead to unconditionally numerically stable discrete systems and provide best approximation properties in terms of the energy norm, $\norm{\cdot}{B}: \UU\longrightarrow [0,\infty]$, of the error, i.e.:
\begin{equation}
\label{eq:energy_norm}
\norm{(u,\qq)}{B} \isdef \supp{(v,\www)\in \VV\setminus \{(0,\mathbf{0})\}} 
     \frac{|B((u,\qq);(v,\www))|}{\norm{(v,\www)}{\VV}}.
\end{equation}

 The discrete fluxes we use, are more regular than is required by the minimal topology we described in the previous section. We apply discrete fluxes that belong to $\SHOO$ rather than $\SHdivO$. Our experience indicates that the numerical solutions we obtain when we use approximations in $\SHdivO$  yield similar accuracy, as long as the domain does not exhibit any re-entrant corners and/or cracks. Convergence is observed in the latter case, but the onset of asymptotic convergence is then generally observed at a higher number  of mesh refinements. Raviart-Thomas discretizations most likely resolve this and will be investigated in an upcoming manuscript.  Currently, using discrete fluxes in $\SCZO$  is certainly less challenging to implement $\SHOO$ partitions on standard meshes. Possibly more importantly, this will allow the use of AVS formulations in commercial simulation software by redefining the user-defined elemental routines.

Let  us now proceed   by deriving the FE discretization of~\eqref{eq:weak_form} by first introducing the   family of invertible maps, $\{\FF:\KK\subset\mathbb{R}^2 \longrightarrow \Omega\}$, such that every $\Kep$ is the image of a master element $\KK$ through one of the mappings $\FF$ (see Figure~\ref{fig:model_problem}). The (conforming) space of  trial functions, $\UUh\subset \UU$, is then defined as:
\begin{equation}
\label{eq:UhPh}
\begin{array}{l}
\ds \UUh \isdef \biggl\{ (\varphi^h,{\pmb \theta}^h)\in \SCZO \times \SCZOvec: \;
(\varphi^h_{|K_m},{\pmb \theta}^h_{|K_m}) = (\hat{\varphi},\hat{\pmb \theta}) \circ \FF , \;\;
                \biggr.
 \\[0.1in]
         \ds \hfill \biggl.  \hat{\varphi} \in P^{p_m}(\KK) \; \wedge\; \hat{\pmb \theta} \in [P^{p_m}(\KK)]^2, \; \;\;
                \forall \Kep
                \biggr\},
\end{array}
\end{equation}
where $p_m$ denotes the local polynomial degree of approximation on $K_m$. We are essentially  following the classical  FE method here and therefore accordingly represent  the FE approximations, $u^h$ and $\qq^h=\{q^h_x,q^h_y\}^T$,  as linear combinations of  trial functions $(e^i(\xx),(E_x^j(\xx), E^k_y(\xx)))\in\UUh$ and corresponding degrees of freedom, $\{u^h_i\in\mathbb{R}, \, i=1,2,\dots, N\}$, $\{q^{h,j}_x\in\mathbb{R}, \, j=1,2,\dots,N\}$ and $\{q^{h,k}_y\in\mathbb{R},\, k=1,2,\dots,N\}$;  i.e., 
\begin{equation} \label{eq:FE_sol}
u^h(\xx) = \summa{i=1}{N} u^h_i \, e^i(\xx), 
\quad q^h_x(\xx) = \summa{j=1}{N} q^{h,j}_x \, E_x^j(\xx),
\quad q^h_y(\xx) = \summa{k=1}{N} q^{h,k}_y \, E_y^k(\xx).
\end{equation}
As mentioned previously, contrary to the trial functions (which are global $C^0$  functions), the test functions are to be piecewise discontinuous and  constructed by invoking the DPG strategy~\cite{Demkowicz4, Demkowicz1, Demkowicz2, Demkowicz3, Demkowicz5, Demkowicz6,NIEMI20132096}. Each of the $3N$ trial functions $e^i(\xx)$, $E_x^j(\xx)$, and $E_y^k(\xx)$, is paired with a vector-valued  test function. Thus,  $e^i(\xx)$ is paired with $(\vvi,\wwwi)\in\VV$,  $E_x^j(\xx)$ with $(\vvj,\wwwj)\in\VV$, and $\E^k_y(\xx)$ with $(\vvk,\wwwk)\in\VV$. Following the DPG philosophy, these pairings are established through  the following  variational problems:
\begin{equation}
\label{eq:test_problems}
\begin{array}{rcll}
\ds \left(\, (r,\zzz);(\vvi,\wwwi) \, \right)_\VV &  \! \! =  \! & B(\,(e^i,\mathbf{0});(r,\zzz) \, ),& \, \forall (r,\zz)\in\VV, \quad i=1,\dots, N,
\\[0.1in]
\ds \left(\, (r,\zzz); (\vvj,\wwwj) \, \right)_\VV & \!  \!  =  \! & B(\, (0,(E_x^j,0));(r,\zzz)\, ), & \, \forall (r,\zz)\in\VV, \quad j=1,\dots, N,
\\[0.1in]
\ds \left(\, (r,\zzz); (\vvk,\wwwk) \, \right)_\VV & \!  \!  =  \! & B(\, (0,(0,E_y^k));(r,\zzz) \,), & \, \forall (r,\zz)\in\VV, \quad k=1,\dots, N,
\end{array}
\end{equation}
where $\left(\, \cdot; \cdot \, \right)_\VV:\; \VV\times\VV \longrightarrow \mathbb{R}$, is the  inner product:
\begin{equation}
\label{eq:inner_prod_V}
\left(\, (r,\zzz); (v,\www) \, \right)_\VV 
\isdef  \summa{\Kep}{}\int_{K_m} \biggl[  h_m^2 \Nabla r_m \cdot \Nabla v_m + r_m\, v_m   + \zzz_m \cdot \www_m\biggr] \dx,
\end{equation}
which induces the norm $\norm{\cdot}{\VV}$, as defined in~\eqref{eq:broken_norms}. The solution of these Riesz representation problems in the test space norm produces the set of test functions that we  use in our variational framework.

\begin{rem}
The variational statements~\eqref{eq:test_problems} are infinite dimensional problems which we approximate numerically. To do so, we compute piecewise discontinuous polynomial approximations $(\vvih,\wwwih)$, $(\vvjh,\wwwjh) $, and $(\vvkh,\wwwkh)$ of $(\vvi,\wwwi)$, $(\vvj,\wwwj) $, and $(\vvk,\wwwk)$, respectively, by applying local polynomial degrees of approximation of order $p_m+\Delta p$.
\end{rem}

\begin{rem}
By applying  functions $(r,\zzz)\in \VV$ in the variational statements of~\eqref{eq:test_problems} that vanish outside a given element $K_m$, the local restriction of the test functions to $K_m$, can easily be computed by solving the following local restrictions of~\eqref{eq:test_problems}:
\begin{equation}
\label{eq:local_test_problems}
\begin{array}{rcll}
\ds \left(\, (r,\zzz);(\vvih,\wwwih) \, \right)_\VVK & = & B_{|K_m}(\,(e^i,\mathbf{0});(r,\zzz) \, ),& \quad \forall (r,\zz)\in\VVK, 
\\[0.1in]
\ds \left(\, (r,\zzz); (\vvjh,\wwwjh) \, \right)_\VVK & = & B_{|K_m}(\, (0,(E_x^j,0));(r,\zzz)\, ), & \quad \forall (r,\zz)\in\VVK, 
\\[0.1in]
\ds \left(\, (r,\zzz); (\vvkh,\wwwkh) \, \right)_\VVK & = & B_{|K_m}(\, (0,(0,E_y^k));(r,\zzz) \,), & \quad \forall (r,\zz)\in\VVK, 
\end{array}
\end{equation}
where $B_{|K_m}(\cdot;\cdot)$ denotes the restriction of  $B(\cdot;\cdot)$ (see~\eqref{eq:B_and_F}) to the element $K_m$ and:
\begin{equation}
\label{eq:local_test_spaces}
\begin{array}{c}
\VVK \isdef \biggl\{ (v,\www)\in \SHOK\times\SLTKvec: \;  \gamma_0^m(v_m)_{|\dKm\cap\GD} =0\biggr\},
\\[0.15in]
\left(\, \cdot; \cdot \, \right)_\VVK:\; \VVK\times\VVK \longrightarrow \mathbb{R}, 
\\[0.15in]
\left(\, (r,\zzz); (v,\www) \, \right)_\VVK 
\isdef  \ds \int_{K_m} \biggl[  h_m^2 \Nabla r \cdot \Nabla v+ r\, v   + \zzz \cdot \www\biggr] \dx.
\end{array}
\end{equation}
If we look at the action of the local restriction of the bilinear form $B(\cdot;\cdot)$ onto functions $(\fr,\vfr)$, that have the same regularity as our FE trial functions (i.e., they belong to $\SCZO\times\SCZOvec$), and test functions $(r,\zz)\in\VV$ that vanish outside $K_m$, we get from~\eqref{eq:B_and_F}:
\begin{equation}\label{eq:local_B}
\begin{array}{c}
B_{|K_m}((\fr,\vfr);(r.\zz)) =
\ds \int_{K_m}\biggl[  \, \left( \vfr_m - \DD\Nabla \fr_m\right) \cdot \zz_m \, 
+  \, \vfr_m \cdot \Nabla r_m \, + \,
 (\bb \cdot \Nabla \fr_m) \, r_m  \biggr] \dx
 \\[0.1in]
 \hfill \ds 
 - \int_{\dKm\setminus \overline{\GD\cup\GN}}  \gamma^m_\nn(\vfr_m) \, \gamma^m_0(r_m) \, \dss
 \end{array}
\end{equation}
Thus, in the computations of the local variational statements of~\eqref{eq:local_test_problems}, the action of $B_{|K_m}(\cdot;\cdot)$ can be applied as shown in~\eqref{eq:local_B}.
\end{rem}

\begin{rem}
Since the action of the bilinear form in the RHS of~\eqref{eq:local_test_problems} is entirely local to the element $K_m$, as given in~\eqref{eq:local_B}, a trial function  only induces a nonzero test function in elements where it has support. Hence,
an additional consequence of~\eqref{eq:local_B}  is that the support of every test function is identical to the support of the corresponding trial  function. 
\end{rem}

At last, the FE discretization of~\eqref{eq:weak_BVP_sum_2}, governing the AVS-FE approximation $(u^h,\qq^h)\in\UUh$ of $(u,\qq)$ can now be introduced as follows:
\begin{equation} \label{eq:discrete_form}
\boxed{
\begin{array}{ll}
\text{Find } &  (u^h,\qq^h) \in \UUh \; \text{ such that:}
\\[0.1in]
 &   B((u^h;\qq^h);(\vvh,\wwwh)) = F((\vvh,\wwwh)), \quad \forall (\vvh,\wwwh)\in \VVh, 
 \end{array}}
\end{equation}
where the finite dimensional subspace of test functions $\VVh\subset\VV$ is spanned by the numerical approximations of the test functions  $\{(\vvih,\wwwih)\}_{i=1}^N$, $\{(\vvjh,\wwwjh)\}_{j=1}^N$, and $\{(\vvkh,\wwwkh)\}_{k=1}^N$, as computed from the Riesz representation problems~\eqref{eq:test_problems} and~\eqref{eq:local_test_problems} by using local polynomial degrees of approximation $p_m+\Delta p$.

Since we essentially apply the DPG methodology~\cite{Demkowicz4, Demkowicz1, Demkowicz2, Demkowicz3, Demkowicz5, Demkowicz6}  in the construction of the space of test functions $\VVh$ via the Riesz representation statements~\eqref{eq:test_problems}, an important consequence is that the FE discretization~\eqref{eq:discrete_form} also inherits the unconditional numerical stability property of the DPG method. Thus, there is no need for any, generally arduous, determination of  problem and mesh dependent stabilization terms to stabilize the numerical scheme, as done in stabilized FE methods such as SUPG, GLS, and VMS. The discrete problem~\eqref{eq:discrete_form} is automatically and unconditionally stable for any choice of the mesh parameters $h_m$ and $p_m$.

\begin{lem}
\label{thm:local_conservation}
The FE discretization~\eqref{eq:discrete_form} is locally conservative.
\end{lem}
We refer to~\cite{CaloRomkes2018} for a detailed  proof of this lemma~\qed

\section{Exemplary Numerical Results}
\label{sec:numerical_verifications}
To conduct numerical studies of our new method, we consider the following simplified form of our model scalar-valued convection diffusion problem~\eqref{eq:conv_diff_BVP} on the unit square domain $\Omega = (0,1)\times(0,1)\subset \mathbb{R}^2$ with homogeneous Dirichlet boundary conditions:
\begin{equation}\label{eq:num_results_homo_ conv_diff}  
\begin{array}{rl}
\ds  - D \bigtriangleup u +
 \bb \cdot  \Nabla u = f, & \quad \text{ in } \Omega, 
 \\[0.1in]
  u = 0, & \quad \text{ on } \partial \Omega,  
 \end{array} 
\end{equation}
where  the coefficient $D\in\SLOinf$ is a scalar-valued isotropic diffusion coefficient.
In the following subsections, we first verify the asymptotic convergence behavior of the newly introduced AVS-FE method. In Section~\ref{sec:num_results_convergence} we analyze a case in which convection is still rather moderate. However, since our main purpose is to investigate the intrinsic stability property of the method, we focus our attention on convection dominated problems in the  subsections that follow. In Section~\ref{sec:numerical_verifications_conv_diff_homo}, we first look at a classical scenario in which all coefficients are homogeneous, i.e., constant, throughout $\Omega$. Next, we consider a scenario of importance to engineering applications. In Section~\ref{sec:num_results_conv_diff_hetero}, the diffusion coefficient is heterogeneous and therefore varies throughout the domain.  Lastly, we briefly investigate the converse situation in Section~\ref{sec:num_results_shock} in which the diffusion is homogeneous, but the convection varies throughout the domain. Particularly, we look at an example in which the variation of the convection coefficient causes the formation of an internal layer.

The purpose of studying these convection-dominated problem is to test the intrinsic (automatic) stability property of the AVS-FE discretizations, which we attained by using the DPG philosophy in the construction of our test functions. We are particularly interested to see if we indeed : 1) obtain automatic stability for any choice of mesh, 2) avoid overly diffused solutions for initial meshes, which is a commonly encountered impediment of LSFEM solution, and 3) avoid solutions with high oscillations at boundary and internal layers that do not tend to attenuate upon mesh refinements, as encountered in classical FE analyses of such problems. 
%
%

\subsection{Asymptotic Convergence Study}
\label{sec:num_results_convergence}

To ascertain the asymptotic convergence rates in terms of the $\SLTO$, $\SHOO$, and $\norm{\cdot}{\UU}$ norms of the error,   we consider a scenario of our model convection diffusion problem~\eqref{eq:num_results_homo_ conv_diff} in which the diffusion coefficient $D=1/Pe$, where we refer to $Pe\in\mathbb{R}^+$ as the Peclet number, and $\bb=\{1,1\}^T$. We choose  $Pe=10$ and the source function $f$  such that the exact theoretical solution is given by:
\begin{equation*} 
u(x,y) = \left[x + \frac{e^{Pe \cdot x}-1}{1-e^{Pe}}\right]\left[y + \frac{e^{Pe \cdot y}-1}{1-e^{Pe}}\right].
\end{equation*}
This solution  exhibits a  boundary layer along the boundaries $x=1$ and $y=1$, but since there is a moderate level of   diffusion (due to the relatively low value of the Peclet number.),  these layers are not sharp. 
%

%
%
\begin{figure}[h]
\subfigure[\label{fig:total_norm_convergence} $\norm{(u,\qq)-(u^h,\qq^h)}{\UUh}$]{\centering
 \includegraphics[width=0.5\textwidth]{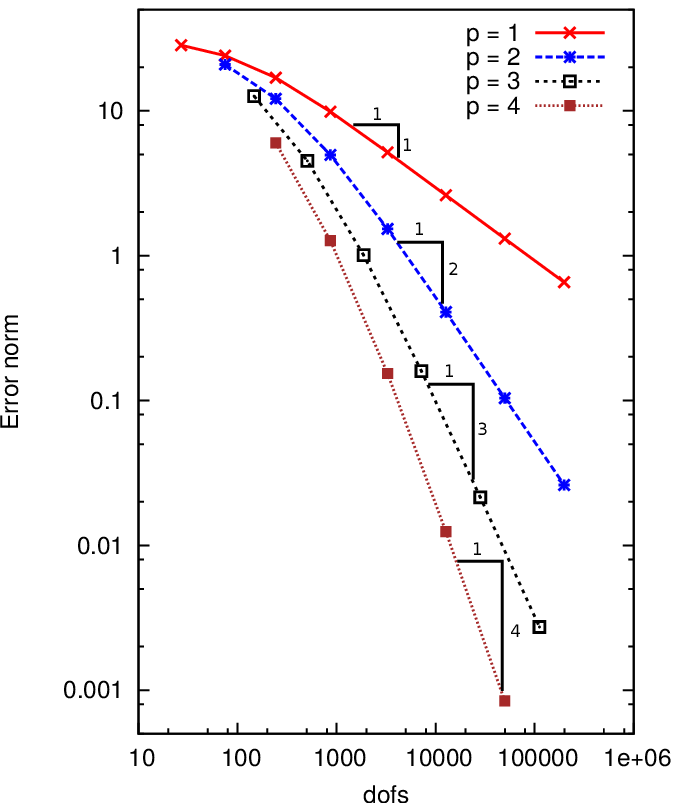}}
\hfill    \subfigure[\label{fig:l2_u_convergence} $\norm{u-u^h}{\SLTO}$]{\centering
 \includegraphics[width=0.4875\textwidth]{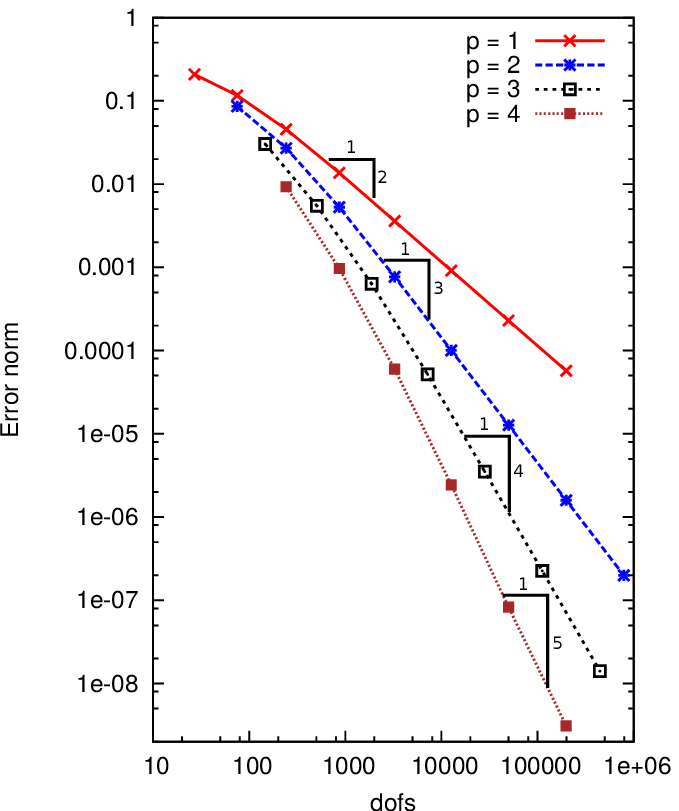}}
 
 \subfigure[\label{fig:h1_u_convergence}$\norm{u-u^h}{\SHOO}$]{\centering
 \includegraphics[width=0.5\textwidth]{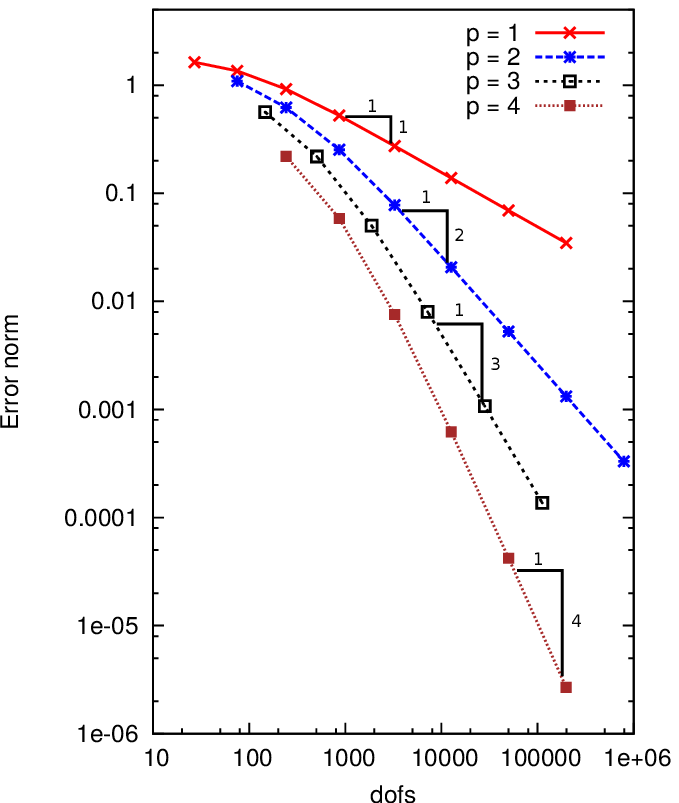}} 
\hfill \subfigure[\label{fig:l2_q_convergence}$\norm{\qq-\qq^h}{\SLTO}$]{\centering
 \includegraphics[width=0.5\textwidth]{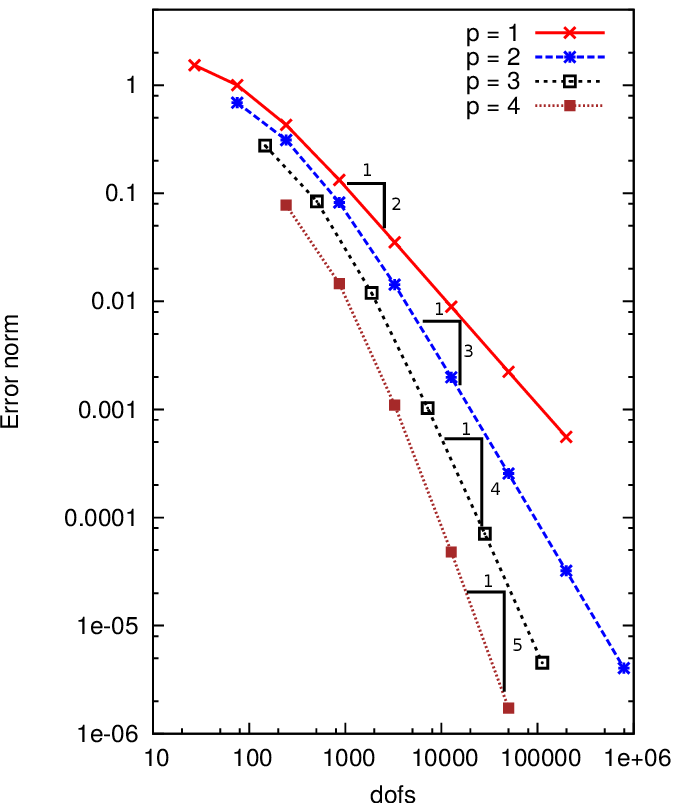}}
  \caption{\label{fig:convergence_results}Error convergence results for uniform $h$-refinements; $\Delta p =0$.}
\end{figure}

%
%
\begin{figure}[t]
\subfigure[\label{fig:q_vs_grad_u_p2} $p=2$]{\centering
 \includegraphics[width=0.5\textwidth]{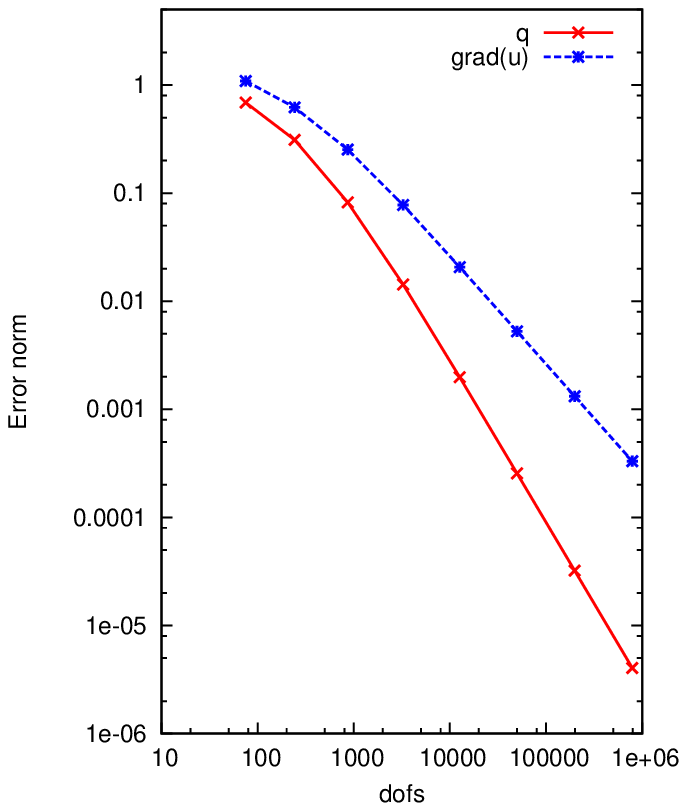}}
\hfill    \subfigure[\label{fig:q_vs_grad_u_p3} $p=3$]{\centering
 \includegraphics[width=0.5\textwidth]{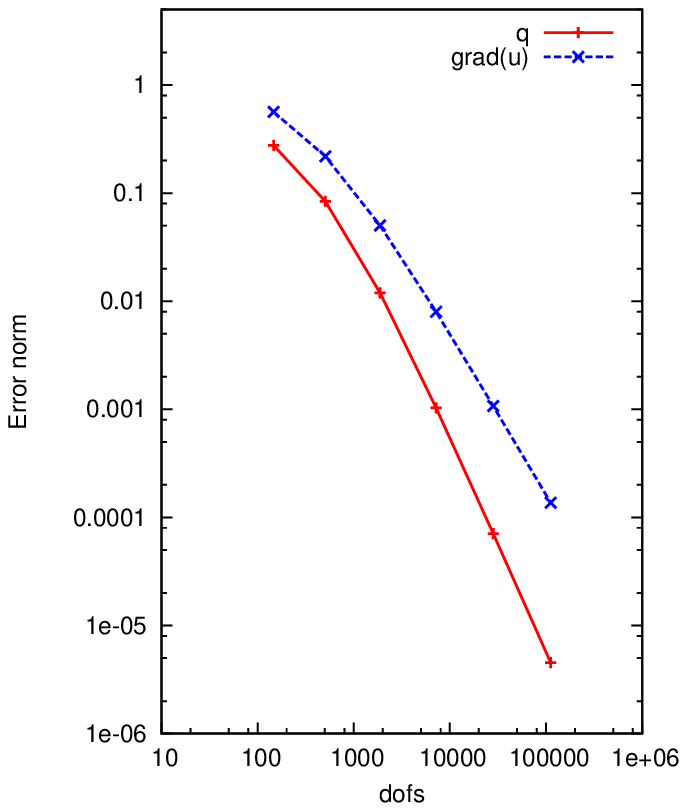}}
   \caption{\label{fig:q_vs_grad_u} Numerical accuracy comparison of $\qq$ versus $\Nabla u$; $\Delta p =0$.}
\end{figure}

In Figure~\ref{fig:convergence_results}, we show   error convergence results for uniform $h$-refinements  in terms of various error norms. For each $h$-refinement study a uniform $p$-level has been applied, ranging from $p=1$ to $p=4$. The test functions have  been computed at the same local polynomial degree of approximation as their corresponding trial functions (i.e., $\Delta p = 0$). The plots in Figures~\ref{fig:l2_u_convergence} and~\ref{fig:h1_u_convergence} clearly show that both the $\SLTO$ and $\SHOO$ norms of the error in the primal variable, $u-u^h$, exhibit optimal convergence rates of order $p+1$ and $p$, respectively. Similarly, the $\SLTO$ norm of the error in the flux, $\qq-\qq^h$, has an optimal convergence rate of $p+1$, as shown in Figure~\ref{fig:l2_q_convergence}. The convergence rates in terms of the error norm $\norm{(u,\qq) - (u^h,\qq^u)}{\UU}$, presented in Figure~\ref{fig:total_norm_convergence}, are also optimal at a rate of $p$. 

These results are representative of extensive convergence studies we have conducted. In all these experiments, the observed asymptotic convergence rates have been optimal. The corresponding a priori estimates of these  convergence rates, and their proofs, are to be presented in~\cite{CaloRomkes2018}.

Lastly, we show a comparison of the $\SLTO$ norm of the error in $\qq-\qq^h$ versus $\Nabla u - \Nabla u^h$ in Figure~\ref{fig:q_vs_grad_u}, for $p=2$ and $p=3$. These results are again representative of extensive numerical experiments, in which consistently a significantly higher accuracy is observed in the prediction of the flux variable versus the gradient of the primal variable.

\subsection{Convection Dominated Diffusion - Homogeneous Coefficients}
\label{sec:numerical_verifications_conv_diff_homo}
%
%

%
%
\begin{figure}[t]
\subfigure[\label{fig:coarse_mesh_homo_conv_diff} Initial $2\times 2$ FE mesh.]{\centering\input{figures/graded_mesh.pstex_t}}
\hfill \subfigure[\label{fig:uh_coarse_homo_conv_diff_p_2} Distribution of $u^h$.]{\centering
 \includegraphics[width=0.6\textwidth]{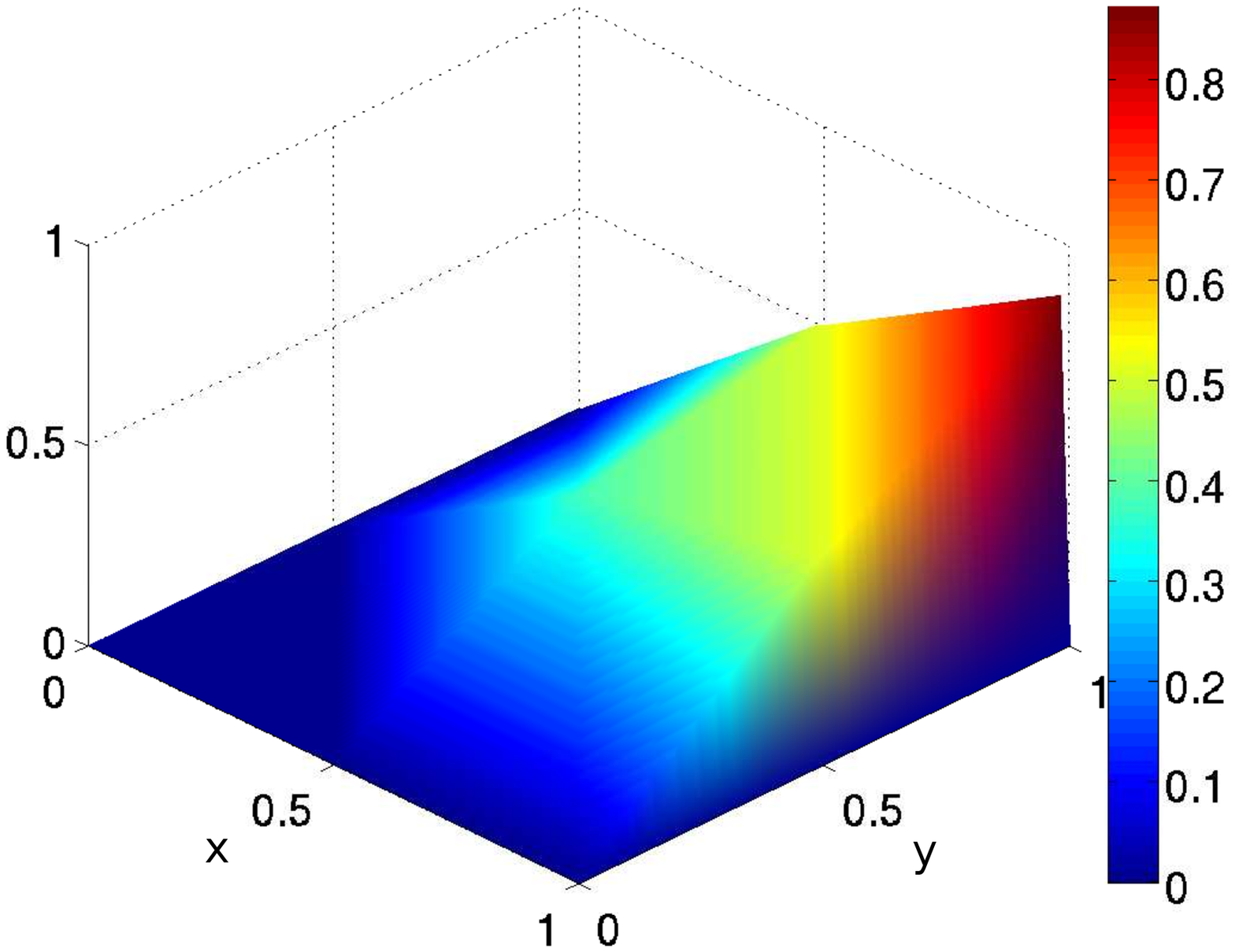} }
  \caption{AVS-FE results  for homogeneous coefficients, a $2\times 2$ graded mesh ($75$ dofs);  $Pe=10^6$,  $p=2$, and $\Delta p =0$.  }
\end{figure}

As mentioned in the introduction,  we are particularly interested in diffusion problems in which convection plays a dominant role.   We start here with the case in which the problem coefficients $D$ and $\bb$ in~\eqref{eq:num_results_homo_ conv_diff} are constant. For our numerical study, we enforce convection  in the diagonal direction, i.e., the convection coefficient $\bb=\{1,1\}^T$. The source function is set at $f(\xx)=1$ and the diffusion coefficient  again at $D=1/Pe$. However,  the Peclet number is now set at a high value of $Pe=10^6$ to ensure the convection term is dominant in~\eqref{eq:num_results_homo_ conv_diff}. With this choice of parameters in place, the distribution of the primal variable  exhibits strong convection in the diagonal direction and a sharp boundary layer of width $1/Pe$ along the boundaries at $x=1$ and $y=1$.

%
%
\begin{figure}[b]
\subfigure[Distribution of $u^h$ throughout $(0,1)\times(0,1)$.]{\centering \includegraphics[width=0.58\textwidth]{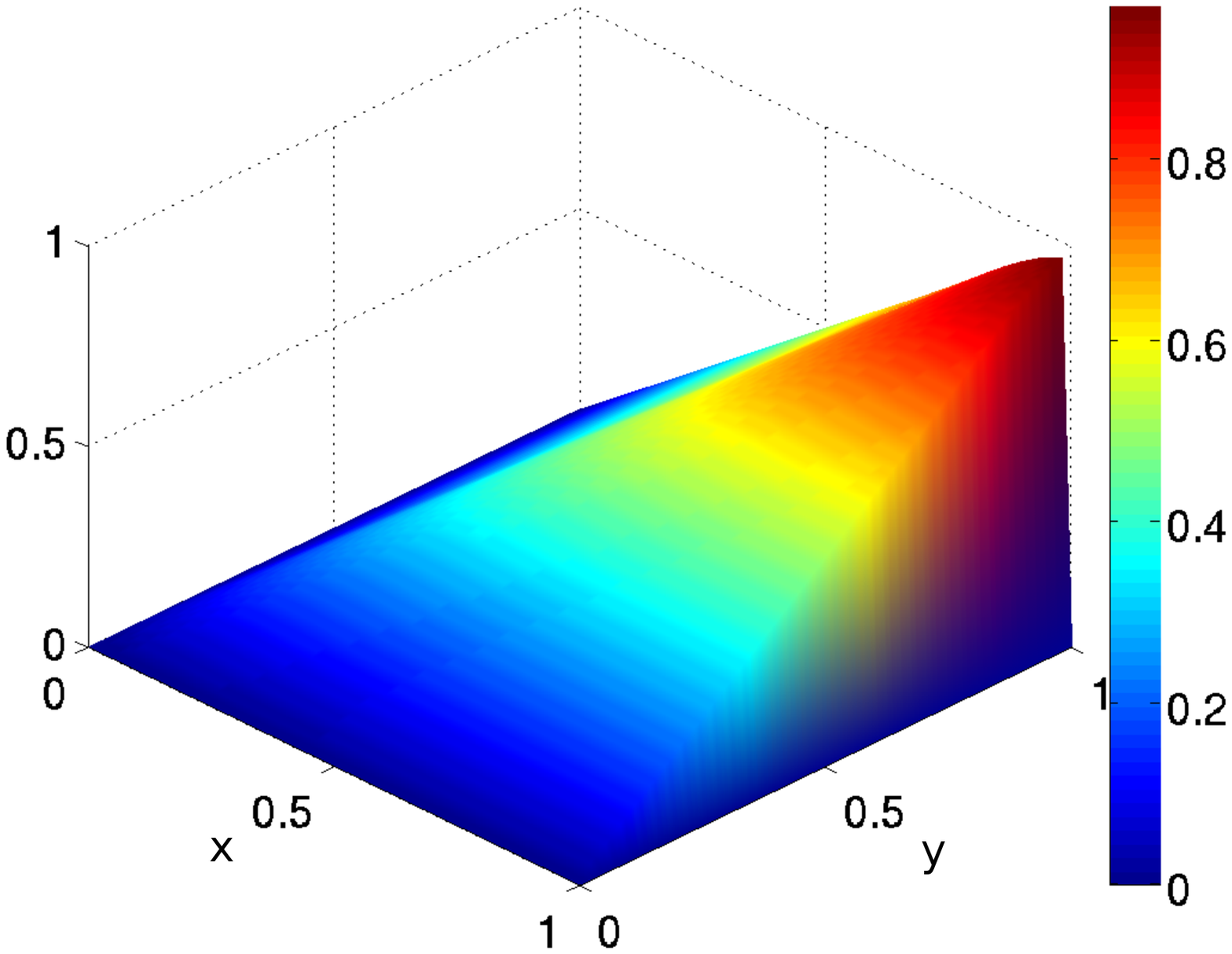}}
\hfill \subfigure[\label{fig:refine_homo_conv_diff_line_plot} Distribution of $u^h$ along  $y-x=0$.]{\centering
 \includegraphics[width=0.4\textwidth]{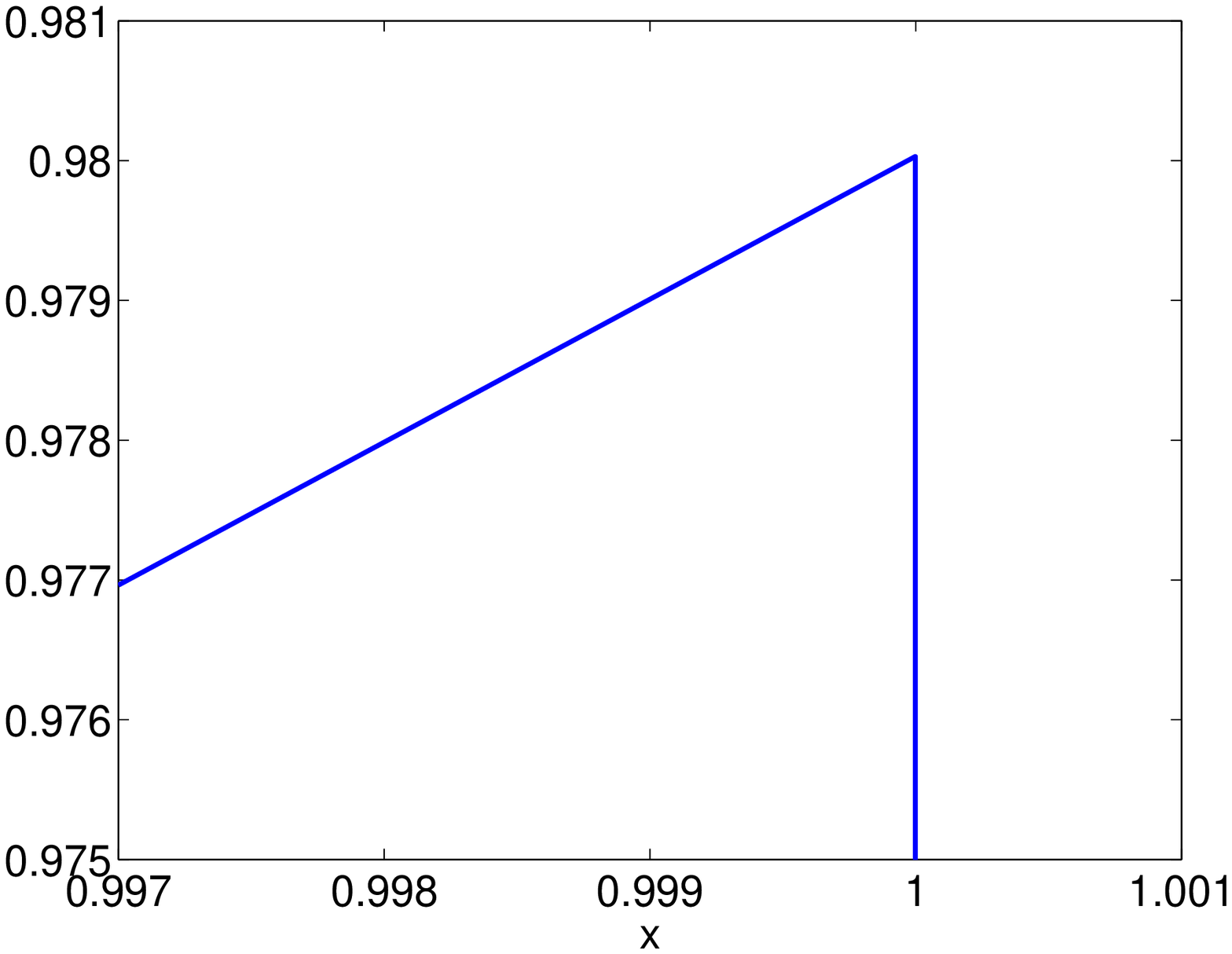}} 
  \caption{\label{fig:refined_homo_conv_diff} AVS-FE results or homogeneous coefficients, a refined mesh  ($12,675$ dofs); $Pe=10^6$, $p=2$, and $\Delta p =0$. }
\end{figure}

For a graded regular mesh of only $2\times 2$ elements, as illustrated in Figure~\ref{fig:coarse_mesh_homo_conv_diff}, and a uniform $p=2$ and $\Delta p =0$, the corresponding   AVS-FE approximation $u^h$   is shown in Figure~\ref{fig:uh_coarse_homo_conv_diff_p_2}.
The FE approximation, at just $75$  dofs, is  stable and does not exhibit  any overly diffused behavior but captures the boundary layer  well. Apparently, applying an identical local polynomial degree of approximation for the test functions in~\eqref{eq:local_test_problems} (i.e. $\Delta p=0$), suffices to   capture the boundary layer with a relatively good accuracy. Results for $\Delta p=1,2,3$, which are not presented here, do not show any significant difference with the results shown here. 

If we apply several additional uniform refinements the solutions remain stable and converge. In Figure~\ref{fig:refined_homo_conv_diff}, results for the AVS-FE approximation are provided for the fifth refinement (i.e., at $12,675$ dofs).  A zoomed-in plot of the distribution of $u^h$ along the diagonal and in the vicinity of the corner at $y=x=1$,  do not show any oscillations, which are commonly observed in solutions obtained via classical FE methods or LSFEM.   The resolution of the boundary layer is not distorted by any oscillations and continuously sharpens as the mesh is refined.

\begin{figure}[t]
\subfigure[\label{fig:unstructured_homo_conv_diff_mesh} Initial unstructured mesh.]{\centering\input{figures/skewed_mesh.pstex_t}}
\hfill \subfigure[\label{fig:unstructured_homo_conv_diff_initial_u} $u^h$ for initial  coarse mesh (27 dofs).]{\centering
 \includegraphics[width=0.45\textwidth]{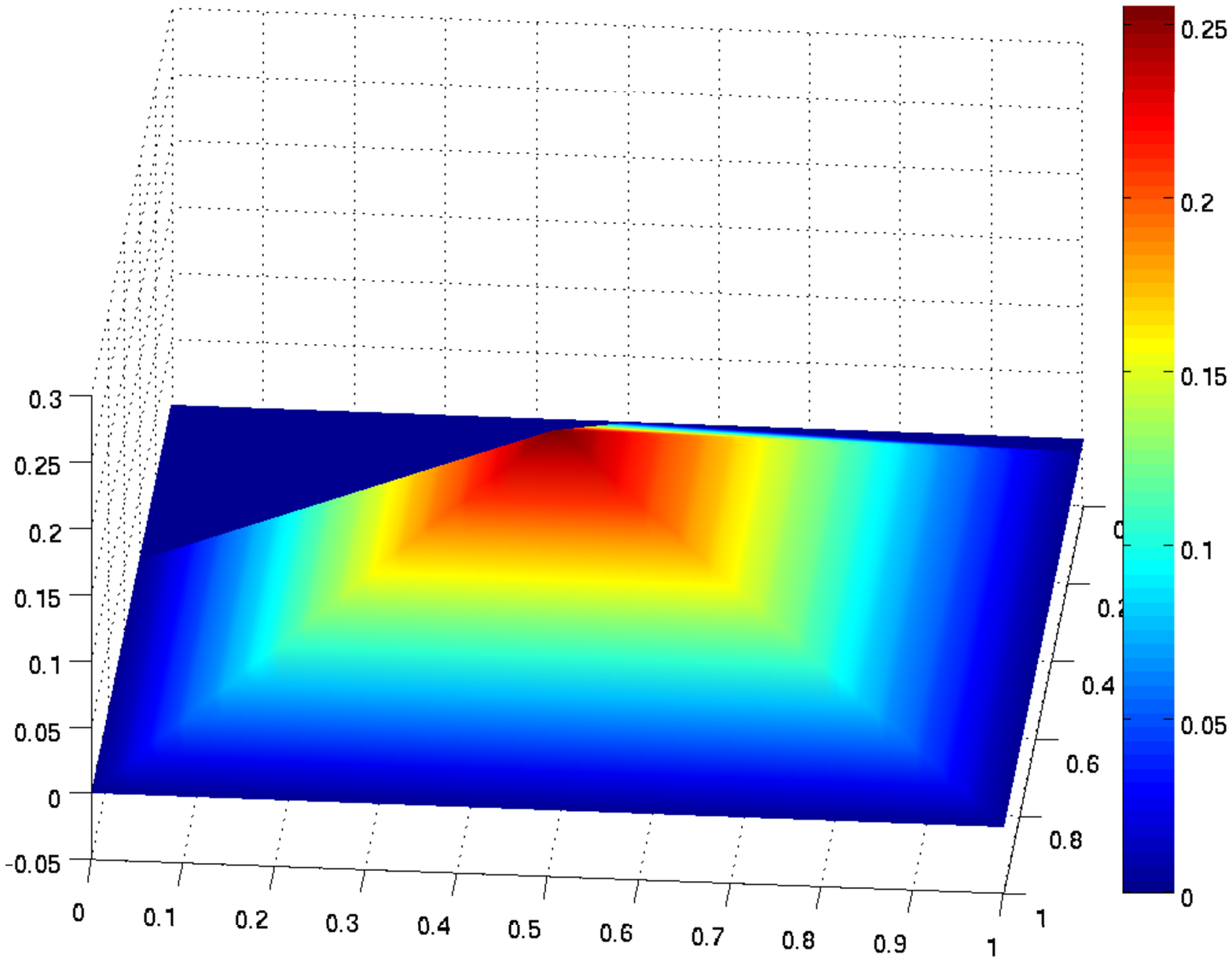}} 
 \\
 \subfigure[\label{fig:unstructured_homo_conv_diff_refined_u} $u^h$ after first  refinement (75dofs).]{\centering
 \includegraphics[width=0.45\textwidth]{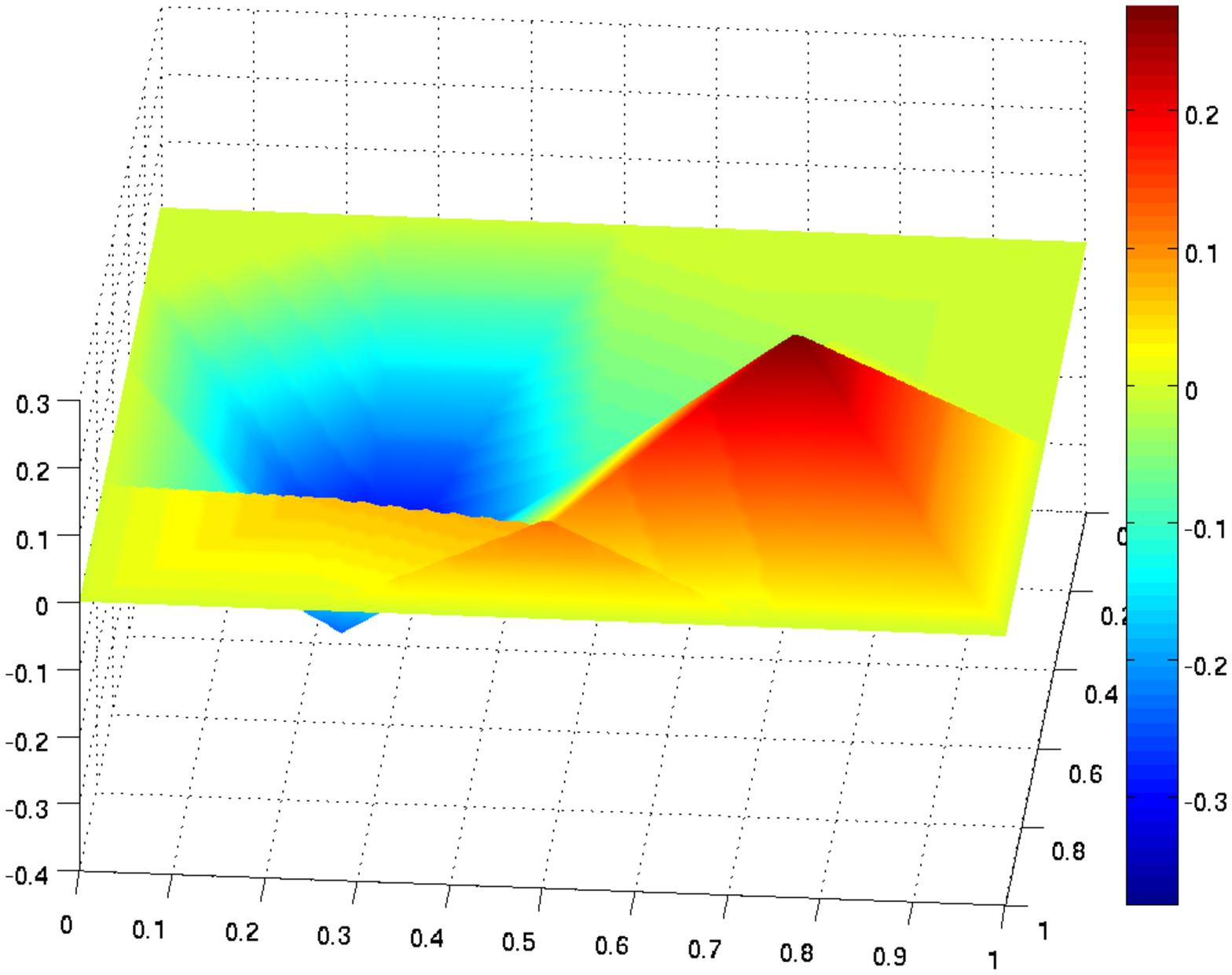}} 
 \hfill \subfigure[\label{fig:unstructured_homo_conv_diff_converged_u} Converged $u^h$.]{\centering
 \includegraphics[width=0.45\textwidth]{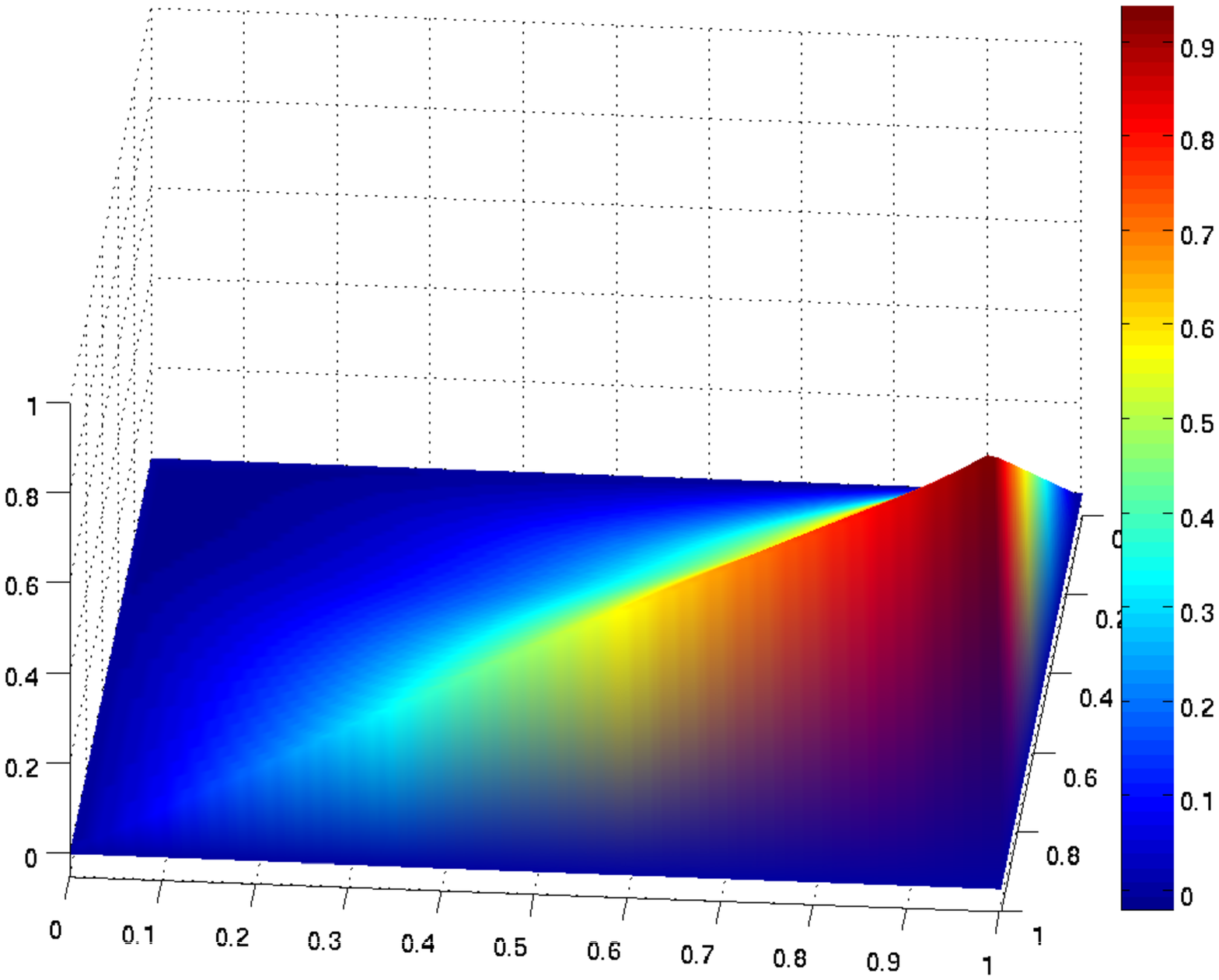}}
  \caption{\label{fig:unstructured_homo_conv_diff} AVS-FE results  for homogeneous coefficients and unstructured  meshes;  $Pe=400$,  $p=2$, and $\Delta p =0$.  }
\end{figure}

To demonstrate that the AVS-FE method also produces sequences of stable numerical solutions for unstructured meshes, we present results in Figure~\ref{fig:unstructured_homo_conv_diff} for $Pe=400$. As depicted in Figure~\ref{fig:unstructured_homo_conv_diff_mesh}, the initial coarse mesh is unstructured and does not resolve the length scale of the boundary layer along $x=1$ and $y=1$. The corresponding numerical solution of $u^h$  is shown in Figure~\ref{fig:unstructured_homo_conv_diff_initial_u} for $p=2$ and has poor numerical accuracy, as is expected for such a coarse mesh. However, the solution is stable and upon applying uniform refinements (see Figure~\ref{fig:unstructured_homo_conv_diff_refined_u} for the first refinement), the solutions indicate the presence of the boundary layer. Hence, any subsequent $hp$-adaptive strategies can then be applied to fully resolve the boundary layer. Since in this work our focus is not on $hp$-adaptivity, we simply apply several uniform $h$-refinements to demonstrate that the solutions do converge for unstructured meshes, as shown in Figure~\ref{fig:unstructured_homo_conv_diff_converged_u}. 
%
%

%
%
\begin{figure}[t]
\subfigure[\label{fig:Pe_checker_board} Diffusion coefficient distribution.]{\centering\input{figures/checkerboard_grey.pstex_t} \hspace{0.5in}}
\hfill \subfigure[\label{fig:coarse_mesh_hetero_conv_diff} Initial  $4\times 4$ graded mesh.]{\centering
 \input{figures/16-elem-graded.pstex_t} }
  \caption{Checker board problem. }
\end{figure}
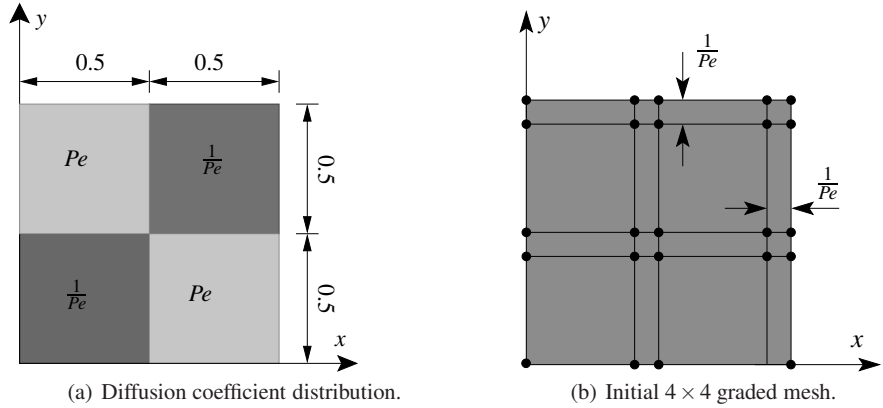

\subsection{Heterogeneous Diffusion}
\label{sec:num_results_conv_diff_hetero}
We continue by looking at a   more challenging case in which the diffusion $D$ is a discontinuous piecewise constant function.  Specifically, $D$ has a value of $Pe$ or $1/Pe$ following a checker board pattern, as depicted in Figure~\ref{fig:Pe_checker_board}. Both the source function and convection coefficient remain unchanged from the experiment conducted in Section~\ref{sec:numerical_verifications_conv_diff_homo}, i.e., $f(\xx)=1$ and $\bb=\{1,1\}^T$. 
By  choosing  a high Peclet number of $10^4$, we essentially establish a  zero solution  in the diffusion dominant quadrants of the domain, while strong convection is observed in the remaining two quadrants. Consequently, in the convective regions, sharp internal layers are formed at the interface with the diffusion dominant quadrants,, with a width of approximately $1/Pe$. Additionally,  sharp boundary layers are present in the convective quadrants along their  boundaries that intersect with the outer boundaries at $x=1$ and $y=1$.

\begin{figure}[t]
\subfigure[\label{fig:uh_coarse_hetero_conv_diff_p_2} Distribution of $u^h$ throughout $(0,1)\times(0,1)$.]{\centering \includegraphics[width=0.6\textwidth]{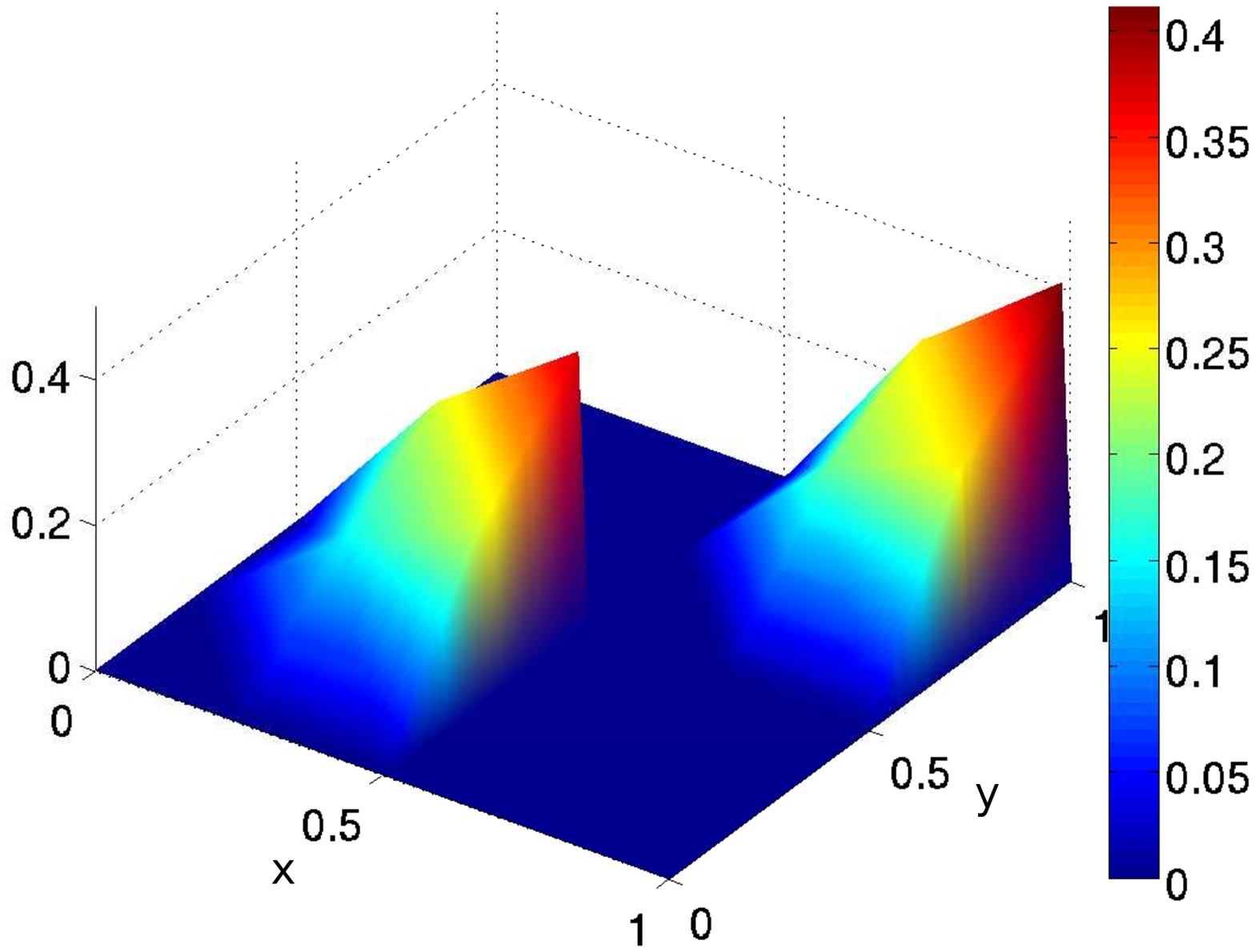}}
\hfill \subfigure[\label{fig:uh_coarse_hetero_conv_diff_p_2_line} $u^h$ along  the diagonal $y-x=0$.]{\centering
 \includegraphics[width=0.4\textwidth]{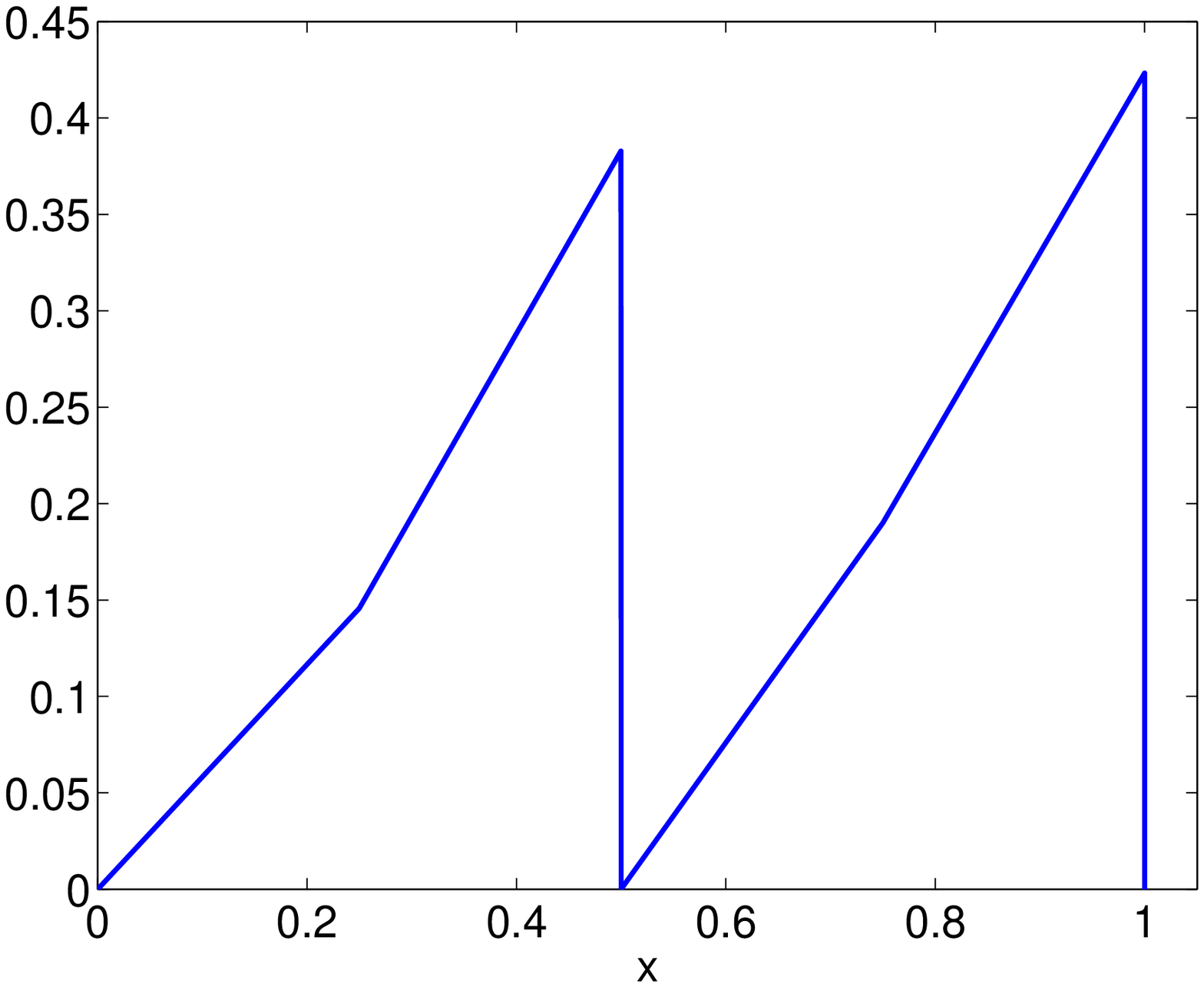}} 
  \caption{\label{fig:coarse_hetero_conv_diff_line_plot} AVS-FE results for heterogeneous diffusion,  $4\times 4$ graded mesh ($243$ dofs), $p=2$, $\Delta p =0$. }
\end{figure}

For a graded regular mesh of only $4\times 4$ elements (see Figure~\ref{fig:coarse_mesh_hetero_conv_diff}),  $p=2$, and $\Delta p =0$, a  contour plot of the  of the distribution of the corresponding AVS-FE solution, $u^h$, throughout the unit square is depicted in   Figure~\ref{fig:uh_coarse_hetero_conv_diff_p_2};  whereas in Figure~\ref{fig:uh_coarse_hetero_conv_diff_p_2_line} its distribution along the diagonal $y-x=0$ is presented. Analogous to the results in Section~\ref{sec:numerical_verifications_conv_diff_homo}, the numerical solution  successfully captures the main features of the solution, i.e., the solution indeed vanishes in the diffusion dominant quadrants, strong convection is seen in the remaining regions, and the sharp internal and boundary layers are adequately captured. It is remarkable that with only $16$ elements, and $243$ dofs, the AVS-FE computation succeeds in resolving these features without any strong oscillations and without the need for any artificial stabilization. Again, using the same polynomial degree of approximation in solving the optimal test functions~\eqref{eq:local_test_problems}, does not appear to inhibit the corresponding AVS-FE computation to resolve the essential solution features.

%
%

\begin{figure}[b]
\subfigure[Distribution of $u^h$ throughout $(0,1)\times(0,1)$.]{\centering \includegraphics[width=0.59\textwidth]{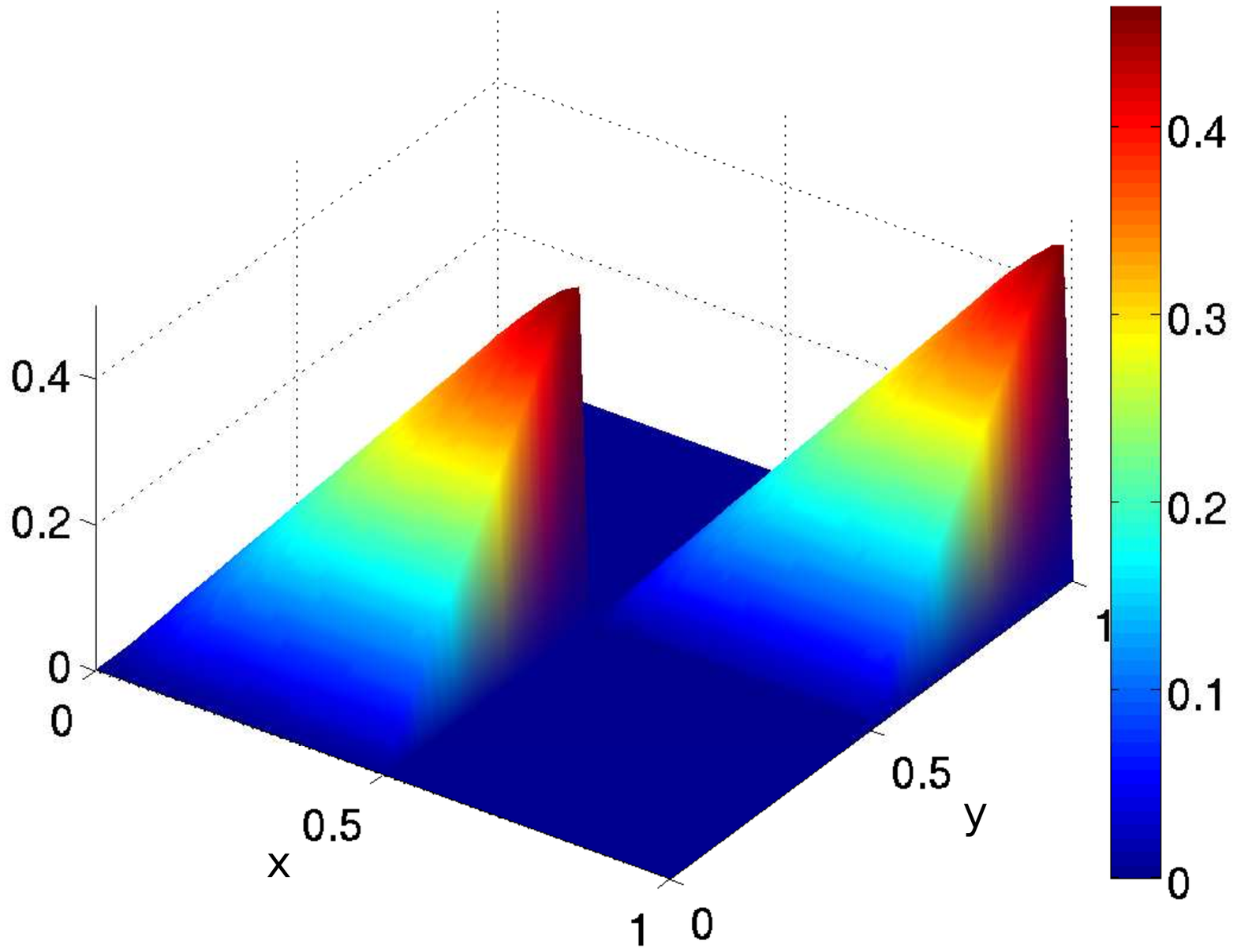}}
\hfill \subfigure[$u^h$ along  the diagonal $y-x=0$.]{\centering
 \includegraphics[width=0.38\textwidth]{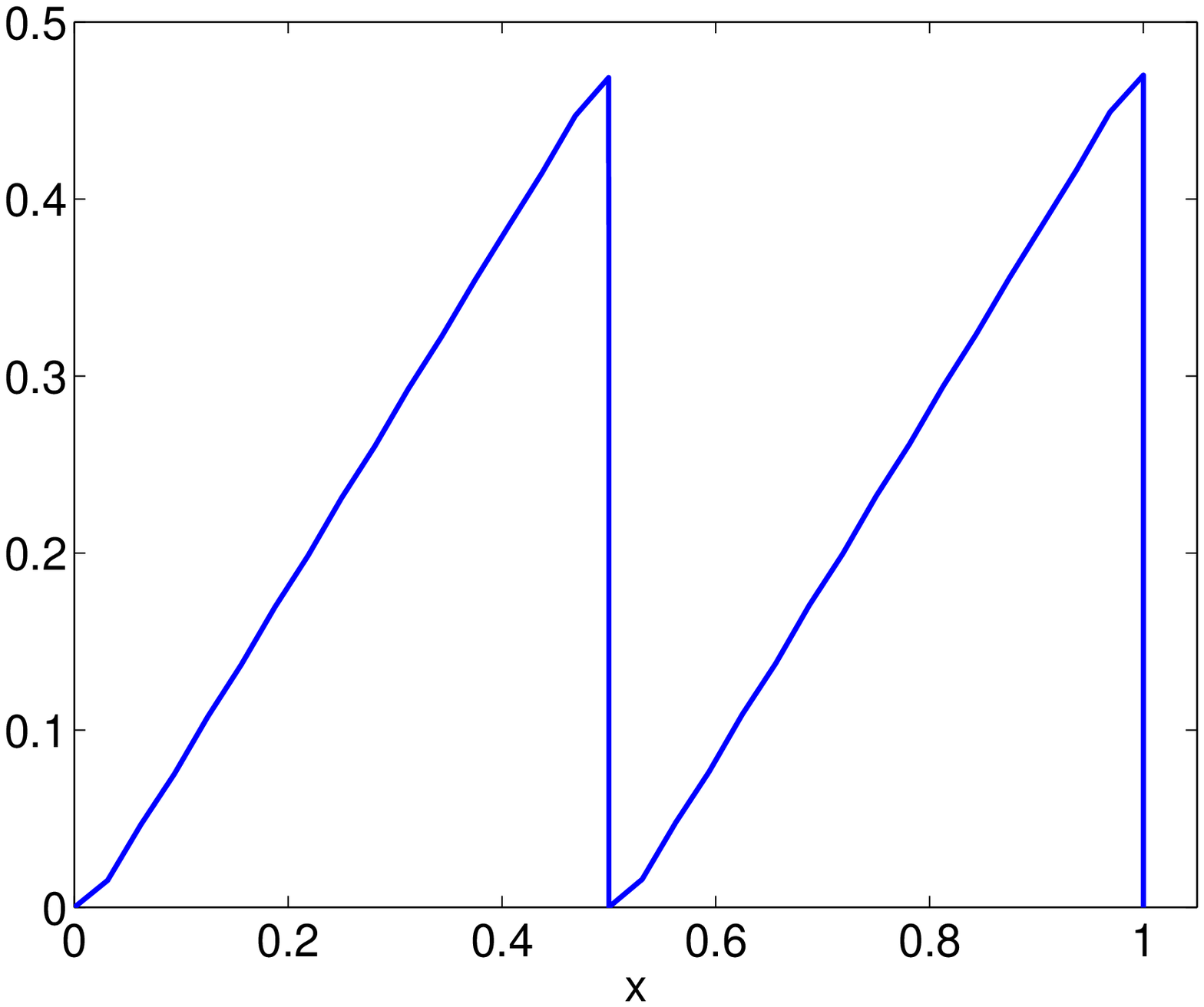}} 
  \caption{\label{fig:refine_hetero_conv_diff_line_plot} AVS-FE results for heterogeneous $D$,  a refined graded mesh ($12,675$ dofs), $p=2$, $\Delta p =0$. }
\end{figure}

Subsequently applying uniform $h$-refinements results in a sequence of numerical solutions, in which the resolution of the internal and boundary continuously improves without inducing any oscillations. Results for the fourth $h$-refinement are given in Figure~\ref{fig:refine_hetero_conv_diff_line_plot}

%
%

\subsection{Non-Constant  Convection}
\label{sec:num_results_shock}
\begin{figure}[t]
\sidecaption
\includegraphics[width=0.6495\textwidth]{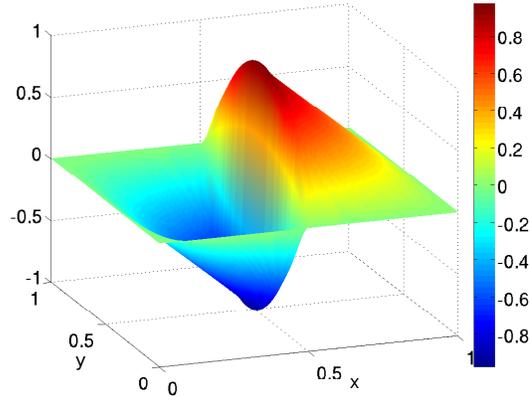}
  \centering \caption{\label{fig:refine_shock} AVS-FE results for non-constant convection,  a refined uniform mesh ($\sim 790,000$ dofs), $Pe=10^9$, $p=1$, and $\Delta p =0$. }
\end{figure}

Lastly, let us now consider a case in which the convection coefficient, $\bb$, rather than the diffusion coefficient, is non-constant, i.e., $\bb = \{ \frac{1}{2}(1-2x),0\}^T$, i.e., we only have convection in the $x$-direction, which varies linearly throughout $\Omega$ and vanishes along the middle line segment $x=\frac{1}{2}$. By choosing the Peclet number at an extremely high level, $Pe=10^9$, we ensure that convection is heavily dominant away from the line segment $x=\frac{1}{2}$. Next, the source function is chosen to be:

\begin{equation}\notag
f(x,y) = \frac{4x-2}{Pe} + y(1-y^2) (8x-4).
\end{equation}
Under these conditions, the solution exhibits a sharp internal layer along the middle line segment $x=\frac{1}{2}$, with a width of the order of $1/Pe$, i.e., $10^{-9}$. Away from the internal layer, or 'shock', the solution is convective. In Figure~\ref{fig:refine_shock}, we present the distribution of $u^h$ for the case in which we started with a $2\times 2$ uniform mesh, $p=1$, and $\Delta p = 0$, and subsequently applied seven uniform $h$-refinements, arriving at a mesh with approximately $790$k dofs. The numerical solutions do not show any oscillatory behavior close to the shock and continuously  provide sharper resolutions of the internal layer as the mesh is refined, while converging to a bounded amplitude. It is striking that the results are automatically stable for a staggering value of a billion for the Peclet number.

\section{Concluding Remarks}
\label{sec:conclusions}

We constructed a variationally stable finite element discretization. This hybrid continuous-discontinuous Petrov-Galerkin method uses solution (trial) functions that are piecewise continuous over the whole domain. That is, these functions correspond to standard finite element partitions. We then use as weight (test) functions a piecewise discontinuous basis. This broken test space allows us to extend the DPG approach to compute optimal test functions automatically and with these to establish numerically stable FE approximations. Important features of this discretization are as follows. The support of each discontinuous test function is identical to its corresponding continuous trial function. The local test-function contribution computed locally on an element by element fashion(i.e. decoupled). This has a linear cost with respect to the problem size and can be thought as an alternative assembly process, where not only inner products, but the functions themselves need to be computed on the fly. Additionally, our experience indicates that the computation of the optimal test functions is achieved with sufficient accuracy by using the same polynomial order of approximation,  p, as that used in the trial function. As in every other DPG formulation, the resulting algebraic system is symmetric and positive definite, allowing us to use simple iterative strategies to compute the numerical solution.  Our future work will include developing variationally stable discretizations based on isogeometric analysis (IGA) both in Galerkin as well as in collocation form. Our preliminary results indicate that these methods are very promising by delivering robust and efficient discretizations exploiting the smoothness of IGA basis functions to deliver intrinsically stable discretizations that are symmetric and positive definite for arbitrary partial differential equations.

We are confident in the impact this methodology will have. Thus, we are partnering with the development communities around FireDrake,  Fenics-HPC and Camellia  as well as our traditional partners who develop PetIGA and PetIGA-MF to release portable parallel implementations of this methodology.

\begin{acknowledgement}
The support of the authors, Albert Romkes and Eirik Valseth, by the NSF CBET  Program, under  NSF Grant  titled \emph{Sustainable System for Mineral Beneficiation}, NSF Grant No. 1805550,  is gratefully acknowledged. This publication was also made possible in part by the CSIRO Professorial Chair in Computational Geoscience at Curtin University and the Deep Earth Imaging Enterprise Future Science Platforms of the Commonwealth Scientific Industrial Research Organisation, CSIRO, of Australia. Additional support was provided by the European Union's Horizon 2020 Research and Innovation Program of the Marie Sk\l{}odowska-Curie grant agreement No. 777778, and the Mega-grant of the Russian Federation Government (N 14.Y26.31.0013).  Additional, support was provided at  Curtin University by The Institute for Geoscience Research (TIGeR) and by the Curtin Institute for Computation. The J. Tinsley Oden Faculty Fellowship Research Program at the Institute for Computational Engineering and Sciences (ICES) of the University of Texas at Austin has partially supported the visits of author, Victor M. Calo, to ICES where he worked closely with Professor Leszek F. Demkowicz.\end{acknowledgement}

 \bibliographystyle{spmpsci}

 \end{document}

%% file: figures/fig-model-with-partition.pstex_t
\begin{picture}(0,0)%
\includegraphics{figures/fig-model-with-partition.pstex}%
\end{picture}%
\setlength{\unitlength}{1973sp}%
\begingroup\makeatletter\ifx\SetFigFont\undefined%
\gdef\SetFigFont#1#2#3#4#5{%
  \reset@font\fontsize{#1}{#2pt}%
  \fontfamily{#3}\fontseries{#4}\fontshape{#5}%
  \selectfont}%
\fi\endgroup%
\begin{picture}(8673,4377)(1561,-4741)
\put(7304,-3775){\makebox(0,0)[lb]{\smash{{\SetFigFont{12}{14.4}{\rmdefault}{\mddefault}{\updefault}{\color[rgb]{0,0,0}$\bfm{F}_n$}%
}}}}
\put(4426,-511){\makebox(0,0)[lb]{\smash{{\SetFigFont{11}{13.2}{\rmdefault}{\mddefault}{\updefault}{\color[rgb]{0,0,0}$\GN$}%
}}}}
\put(1576,-3661){\makebox(0,0)[lb]{\smash{{\SetFigFont{11}{13.2}{\rmdefault}{\mddefault}{\updefault}{\color[rgb]{0,0,0}$\GD$}%
}}}}
\put(3976,-2461){\makebox(0,0)[lb]{\smash{{\SetFigFont{12}{14.4}{\rmdefault}{\mddefault}{\updefault}{\color[rgb]{0,0,0}$K_m$}%
}}}}
\put(9901,-1861){\makebox(0,0)[lb]{\smash{{\SetFigFont{12}{14.4}{\rmdefault}{\mddefault}{\updefault}{\color[rgb]{0,0,0}$\xi$}%
}}}}
\put(9001,-1036){\makebox(0,0)[lb]{\smash{{\SetFigFont{12}{14.4}{\rmdefault}{\mddefault}{\updefault}{\color[rgb]{0,0,0}$\eta$}%
}}}}
\put(5401,-2911){\makebox(0,0)[lb]{\smash{{\SetFigFont{12}{14.4}{\rmdefault}{\mddefault}{\updefault}{\color[rgb]{0,0,0}$K_n$}%
}}}}
\put(6601,-792){\makebox(0,0)[lb]{\smash{{\SetFigFont{12}{14.4}{\rmdefault}{\mddefault}{\updefault}{\color[rgb]{0,0,0}$\bfm{F}_m$}%
}}}}
\put(5385,-4450){\makebox(0,0)[lb]{\smash{{\SetFigFont{12}{14.4}{\rmdefault}{\mddefault}{\updefault}{\color[rgb]{0,0,0}$\nn$}%
}}}}
\put(2543,-934){\makebox(0,0)[lb]{\smash{{\SetFigFont{12}{14.4}{\rmdefault}{\mddefault}{\updefault}{\color[rgb]{0,0,0}$\Omega$}%
}}}}
\put(9020,-1703){\makebox(0,0)[lb]{\smash{{\SetFigFont{12}{14.4}{\rmdefault}{\mddefault}{\updefault}{\color[rgb]{0,0,0}$\KK$}%
}}}}
\end{picture}%

%% file: figures/graded_mesh.pstex_t
\begin{picture}(0,0)%
\includegraphics{figures/graded_mesh.pstex}%
\end{picture}%
\setlength{\unitlength}{3315sp}%
\begingroup\makeatletter\ifx\SetFigFont\undefined%
\gdef\SetFigFont#1#2#3#4#5{%
  \reset@font\fontsize{#1}{#2pt}%
  \fontfamily{#3}\fontseries{#4}\fontshape{#5}%
  \selectfont}%
\fi\endgroup%
\begin{picture}(2575,2529)(1263,-5323)
\put(2296,-3031){\makebox(0,0)[lb]{\smash{{\SetFigFont{10}{12.0}{\rmdefault}{\mddefault}{\updefault}{\color[rgb]{0,0,0}$\frac{1}{Pe}$}%
}}}}
\put(3511,-3976){\makebox(0,0)[lb]{\smash{{\SetFigFont{10}{12.0}{\rmdefault}{\mddefault}{\updefault}{\color[rgb]{0,0,0}$\frac{1}{Pe}$}%
}}}}
\put(3061,-3076){\makebox(0,0)[lb]{\smash{{\SetFigFont{10}{12.0}{\rmdefault}{\mddefault}{\updefault}{\color[rgb]{0,0,0}$\Ph$}%
}}}}
\put(1396,-2986){\makebox(0,0)[lb]{\smash{{\SetFigFont{10}{12.0}{\rmdefault}{\mddefault}{\updefault}{\color[rgb]{0,0,0}$y$}%
}}}}
\put(3511,-5191){\makebox(0,0)[lb]{\smash{{\SetFigFont{10}{12.0}{\rmdefault}{\mddefault}{\updefault}{\color[rgb]{0,0,0}$x$}%
}}}}
\end{picture}%

%% file: figures/skewed_mesh.pstex_t
\begin{picture}(0,0)%
\includegraphics{figures/skewed_mesh.pstex}%
\end{picture}%
\setlength{\unitlength}{3315sp}%
\begingroup\makeatletter\ifx\SetFigFont\undefined%
\gdef\SetFigFont#1#2#3#4#5{%
  \reset@font\fontsize{#1}{#2pt}%
  \fontfamily{#3}\fontseries{#4}\fontshape{#5}%
  \selectfont}%
\fi\endgroup%
\begin{picture}(2585,2585)(1253,-5334)
\put(2066,-5161){\makebox(0,0)[lb]{\smash{{\SetFigFont{10}{12.0}{\rmdefault}{\mddefault}{\updefault}{\color[rgb]{0,0,0}$45^{\circ}$}%
}}}}
\put(1861,-4152){\makebox(0,0)[lb]{\smash{{\SetFigFont{10}{12.0}{\rmdefault}{\mddefault}{\updefault}{\color[rgb]{0,0,0}$30^{\circ}$}%
}}}}
\put(1396,-2986){\makebox(0,0)[lb]{\smash{{\SetFigFont{10}{12.0}{\rmdefault}{\mddefault}{\updefault}{\color[rgb]{0,0,0}$y$}%
}}}}
\put(3511,-5191){\makebox(0,0)[lb]{\smash{{\SetFigFont{10}{12.0}{\rmdefault}{\mddefault}{\updefault}{\color[rgb]{0,0,0}$x$}%
}}}}
\end{picture}%

%% file: figures/checkerboard_grey.pstex_t
\begin{picture}(0,0)%
\includegraphics{figures/checkerboard_grey.pstex}%
\end{picture}%
\setlength{\unitlength}{3108sp}%
\begingroup\makeatletter\ifx\SetFigFont\undefined%
\gdef\SetFigFont#1#2#3#4#5{%
  \reset@font\fontsize{#1}{#2pt}%
  \fontfamily{#3}\fontseries{#4}\fontshape{#5}%
  \selectfont}%
\fi\endgroup%
\begin{picture}(2799,2964)(1219,-5353)
\put(3691,-3616){\rotatebox{270.0}{\makebox(0,0)[lb]{\smash{{\SetFigFont{9}{10.8}{\rmdefault}{\mddefault}{\updefault}{\color[rgb]{0,0,0}0.5}%
}}}}}
\put(1666,-3706){\makebox(0,0)[lb]{\smash{{\SetFigFont{9}{10.8}{\rmdefault}{\mddefault}{\updefault}{\color[rgb]{0,0,0}$Pe$}%
}}}}
\put(2656,-4786){\makebox(0,0)[lb]{\smash{{\SetFigFont{9}{10.8}{\rmdefault}{\mddefault}{\updefault}{\color[rgb]{0,0,0}$Pe$}%
}}}}
\put(1666,-4786){\makebox(0,0)[lb]{\smash{{\SetFigFont{9}{10.8}{\rmdefault}{\mddefault}{\updefault}{\color[rgb]{0,0,0}$\frac{1}{Pe}$}%
}}}}
\put(1756,-2941){\makebox(0,0)[lb]{\smash{{\SetFigFont{9}{10.8}{\rmdefault}{\mddefault}{\updefault}{\color[rgb]{0,0,0}0.5}%
}}}}
\put(2701,-2941){\makebox(0,0)[lb]{\smash{{\SetFigFont{9}{10.8}{\rmdefault}{\mddefault}{\updefault}{\color[rgb]{0,0,0}0.5}%
}}}}
\put(1441,-2581){\makebox(0,0)[lb]{\smash{{\SetFigFont{9}{10.8}{\rmdefault}{\mddefault}{\updefault}{\color[rgb]{0,0,0}$y$}%
}}}}
\put(3826,-5146){\makebox(0,0)[lb]{\smash{{\SetFigFont{9}{10.8}{\rmdefault}{\mddefault}{\updefault}{\color[rgb]{0,0,0}$x$}%
}}}}
\put(2746,-3706){\makebox(0,0)[lb]{\smash{{\SetFigFont{9}{10.8}{\rmdefault}{\mddefault}{\updefault}{\color[rgb]{0,0,0}$\frac{1}{Pe}$}%
}}}}
\put(3691,-4651){\rotatebox{270.0}{\makebox(0,0)[lb]{\smash{{\SetFigFont{9}{10.8}{\rmdefault}{\mddefault}{\updefault}{\color[rgb]{0,0,0}0.5}%
}}}}}
\end{picture}%

%% file: figures/16-elem-graded.pstex_t
\begin{picture}(0,0)%
\includegraphics{figures/16-elem-graded.pstex}%
\end{picture}%
\setlength{\unitlength}{3315sp}%
\begingroup\makeatletter\ifx\SetFigFont\undefined%
\gdef\SetFigFont#1#2#3#4#5{%
  \reset@font\fontsize{#1}{#2pt}%
  \fontfamily{#3}\fontseries{#4}\fontshape{#5}%
  \selectfont}%
\fi\endgroup%
\begin{picture}(2724,2760)(1249,-5338)
\put(2566,-2986){\makebox(0,0)[lb]{\smash{{\SetFigFont{10}{12.0}{\rmdefault}{\mddefault}{\updefault}{\color[rgb]{0,0,0}$\frac{1}{Pe}$}%
}}}}
\put(3466,-3976){\makebox(0,0)[lb]{\smash{{\SetFigFont{10}{12.0}{\rmdefault}{\mddefault}{\updefault}{\color[rgb]{0,0,0}$\frac{1}{Pe}$}%
}}}}
\put(3736,-5146){\makebox(0,0)[lb]{\smash{{\SetFigFont{10}{12.0}{\rmdefault}{\mddefault}{\updefault}{\color[rgb]{0,0,0}$x$}%
}}}}
\put(1396,-2761){\makebox(0,0)[lb]{\smash{{\SetFigFont{10}{12.0}{\rmdefault}{\mddefault}{\updefault}{\color[rgb]{0,0,0}$y$}%
}}}}
\end{picture}%